\documentclass[a4paper,12pt]{article}
\usepackage{times}
\usepackage{amsfonts}
\usepackage{amssymb}
\usepackage{amsmath}
\usepackage{amsthm}
\usepackage[frenchb]{babel}

\begin{document}

\newcommand{\eps}{\varepsilon}

\newtheorem{thm}{Th\'eor\`eme}
\newtheorem*{thm*}{Th\'eor\`eme}
\newtheorem{de}{D\'efinition}[section]
\newtheorem{prop}{Proposition}[section]
\newtheorem{lem}[prop]{Lemme}
\newtheorem{cor}[prop]{Corollaire}

\title{Cibles r\'etr\'ecissantes de rayon $n^{-\frac{1}{d}}$ : propri\'et\'e du logarithme}
\author{Benjamin Mussat}
\date{\today}
\maketitle

\begin{abstract}
On dit qu'une translation sur le tore de dimension $d$ poss\`ede la propri\'et\'e du logarithme si la propri\'et\'e des cibles r\'etr\'ecissantes est v\'erifi\'ee dans le cas des boules de rayon $n^{-\frac{1}{d}}$. En dimension 1, toute rotation irrationnelle poss\`ede la propri\'et\'e du logarithme. En dimension sup\'erieure, nous donnons des crit\`eres permettant de d\'eterminer si cette propri\'et\'e est v\'erifi\'ee ou non. Ces crit\`eres reposent sur une notion de type diophantien diff\'erente de la notion standard.

A l'aide d'une construction en dimension 2 de vecteurs dont nous con\-tr\^ol\-ons les types diophantiens, nous obtenons des contre-exemples \`a la propri\'et\'e du logarithme pour lesquels les vecteurs de translation sont diophantiens d'exposants arbitrairement petits et des exemples poss\'edant la propri\'et\'e du logarithme pour lesquels les vecteurs sont Liouville.

\end{abstract}

%1
\section{Introduction}

Soit $(M,\mathcal{B},\mu, T)$ un syst\`eme dynamique ergodique probabilis\'e inversible, o\`u $M$ est un espace compact m\'etrique et $\mathcal{B}$ l'ensemble de ses bor\'eliens.
%1.1
\subsection{Cibles r\'etr\'ecisssantes}
On suit les d\'efinitions de \cite{Fayad1}.
\begin{de}

Une suite d'ensembles mesurables $\mathcal{A} = (A_n)_{n \in \mathbb{N}}$ est appel\'ee \emph{une suite de Borel-Cantelli (BC) pour T} si, pour presque tout $x$ dans $M$, pour une infinit\'e de $n$, $ T^n(x)$ appartient \`a $A_n$ autrement dit si
$$ \mu(\limsup_{n \to \infty} T^{-n} A_n)=1.$$
\end{de}

\noindent
D'apr\`es le lemme de Borel-Cantelli, il est n\'ecessaire que $$\sum \mu(A_n) = \infty.$$

\begin{lem}
\label{borelcantelli}
Une suite $(A_n)_{n \in \mathbb{N}}$  d\'ecroissante est une suite de Borel-Cantelli si $$ \mu(\limsup_{n \to \infty} T^{-n} A_n) > 0.$$
\end{lem}

On note $B= \limsup(T^{-n} A_n)$. La suite $(A_n)_{n \in \mathbb{N}}$ est d\'ecroissante, donc quel que soit $n$ entier positif on a $T^{-(n+1)} A_{n+1} \subset T^{-1} T^{-n} A_n$, d'o\`u $B \subset T^{-1} B$.
%$$\cup_{n \ge n_0+1} T^{-n} A_n \subset T^{-1}( \cup_{n \ge n_0+1} T^{-n} A_n),$$
%On a donc,
%$$\limsup(T^{-n}A_n) \subset T^{-1} \limsup(T^{-n}A_n).$$
Les $T^{-i} (T^{-1}B \setminus B) $  sont deux \`a deux disjoints et de m\^eme mesure, d'o\`u $\mu(T^{-1}B \setminus B)=0$. L'ensemble $B$ est invariant par $T$ \`a un ensemble de mesure nulle pr\`es.
% $ T^{-1} \limsup A_n \setminus \limsup A_n$ est alors de mesure nulle

Le syst\`eme \'etant suppos\'e ergodique, $B$ est  de mesure nulle ou \'egale \`a un.
\medskip

On dit que le syst\`eme $(M,\mathcal{B},\mu, T)$ a la \emph{propri\'et\'e des cibles r\'etr\'ecissantes } si pour tout $x_0 \in M$, toute suite de boules de centre $x_0$ dont la s\'erie des mesures diverge est BC pour $T$ et qu'il a la \emph{propri\'et\'e des cibles r\'etr\'ecissantes monotone} si pour tout $x_0 \in M$, toute suite d\'ecroissante de boules de centre $x_0$ dont la s\'erie des mesures diverge est BC pour $T$.
%
%1.2
\subsection{Translation sur le tore $\mathbb{T}^d$}
On s'int\'eressera dans toute la suite au syst\`eme ergodique $(\mathbb{T}^d,\mu,T_\theta)$, o\`u $\mathbb{T}^d$ est le tore de dimension $d$ ($d \ge 1$), $\mu$ est la mesure de Lebesgue et $T_\theta$ est la translation par un vecteur $\theta$ dont les coordonn\'ees sont rationnellement ind\'ependantes modulo un.
Un th\'eor\`eme prouv\'e par Kurzweil en 1955 (\cite{Kurzweil}) et red\'ecouvert par Fayad (\cite{Fayad1}) donne d'une part qu'aucune translation n'a la propri\'et\'e des cibles r\'etr\'ecissantes et d'autre part qu'elle poss\`ede la propri\'et\'e des cibles r\'etr\'ecissantes monotone si et seulement si son vecteur est de type constant.
\medskip

Il est naturel de consid\'erer  le cas limite des boules de rayon $n^{-{\frac{1}{d}}} $, puisque nous nous int\'eressons \`a des cibles dont la s\'erie des mesures diverge. S'agissant d'une translation, le choix des $x_0$ est indiff\'erent, on se restreint \`a des boules de centre $0$.

%En dimension $1$, il n'est pas difficile de montrer que la suite de boules de rayon $n^{-1}$ est de Borel-Cantelli pour tout $T_{\theta}$ avec $\theta$ irrationnel, (ce que l'on retrouvera avec nos th\'eor\`eme $1$ et $2$).

\begin{de}
On dit qu'une translation $T_{\theta}$ poss\`ede \emph{la propri\'et\'e du logarithme} si la suite
$( B( x_0 , n^{-{\frac{1}{d}}} ) )_{n \in \mathbb{N}^*}$ est de Borel-Cantelli pour $T_{\theta}$.
\end{de}
Nous utilisons dans la suite la distance sur le tore d\'efinie ci-dessous, le lecteur se convaincra facilement que nos r\'esultats restent vrais pour les distances \'equivalentes.

%
%En effet, si la suite $\left( B( 0 , n^{-{\frac{1}{d}}} ) \right)_{n \in \mathbb{N}^*}$ est de Borel-Cantelli, alors les uistes de boules de m\^eme rayon, centr\'ees en tout point du tore sont de Borel-Cantelli.
\medskip

Un probl\`eme li\'e est celui de ``la loi du logarithme''. On en donne une d\'efinition et on explicite ses liens avec la la propri\'et\'e du logarithme dans la partie 1.5. Un article de Galatolo et Peterlongo (\cite{Galatolo1}) revient sur les liens entre la loi du logarithme et les diff\'erents probl\`emes des temps d'approche d'un point $x_0$ par les orbites d'un autre point.

 .
%1.3
\subsection{Notations et d\'efinitions}
Pour $x$ r\'eel, on notera par $[x]$ la partie enti\`ere de $x$. Pour $x = (x_1,...,x_d) \in \mathbb{R}^d$, on note $|x|=\max_{i=1..d} |x_i|$. La notation $\| .\|$ sera utilis\'e pour la distance \`a l'entier le plus proche dans $\mathbb{R}$ ou au point de $\mathbb{Z}^d$ le plus proche dans $\mathbb{R}^d$ (pour la norme $|.|$). On utilisera la distance donn\'ee par $ \| x- y \| $ entre deux points $x$ et $y$ du tore $\mathbb{T}^d$.

Dans toute la suite, on consid\`ere $\theta= (\theta_1,...,\theta_d)$ un vecteur de $\mathbb{R}^d$ \`a coordonn\'ees rationnellement ind\'ependantes modulo un. On utilise deux approximations de $\theta$, d'une part l'approximation lin\'eaire, o\`u pour $\Delta=(s_1,...,s_d) \in \mathbb{Z}^d$, on consid\`ere
$$ \| \langle \Delta,\theta \rangle \|=\inf_{p \in \mathbb{Z}} \left | s_1\theta_1+...+s_d\theta_d - p \right | ;$$
et d'autre part l'approximation simultan\'ee, o\`u pour $q \in \mathbb{Z}$, on consid\`ere
$$\| q\theta\| =\max_{1 \le i \le d} \ \inf_{p \in \mathbb{Z}} \left | q \theta_{i} -p \right |.$$

Evidemment en dimension $1$, les deux approximations sont confondues.
\medskip

On dit qu'un vecteur $ \Delta$ de  $\mathbb{Z}^d $, non nul, est une \textit{meilleure approximation lin\'eaire} de $\theta \in \mathbb{R}^d$ si  pour tout $\Delta' \in \mathbb{Z}^d$ v\'erifiant $ 0<|\Delta' | <  | \Delta |$, on a

$$ \| \langle \Delta,\theta \rangle \| < \| \langle \Delta',\theta \rangle \|. $$

Notons que $ \| \langle \Delta,\theta \rangle \| = \| \langle \Delta',\theta \rangle \|$ implique $\Delta = \pm \Delta'$, il existe donc une suite $(\Delta_n)_{n \in \mathbb{N}}$, rang\'ee par normes strictement croissantes, compos\'ee, au signe pr\`es, de toutes les meilleures approximation lin\'eaires. On l'appelle suite des meilleures approximations lin\'eaires.

On dit qu'un entier $q$, strictement positif, est \textit{une meilleure approximation simultan\'ee} de $\theta \in \mathbb{R}^d$ si que quel que soit $q'$ entier tel que $0< q'  < q$, on a
$$ \| q\theta\| < \| q'\theta \|.$$
On appelle suite des meilleures approximations simultan\'ees la suite strictement croissante, not\'ee $(q_n) _{n \in \mathbb{N}}$, qui est compos\'ee de toutes les meilleures approximations simultan\'ees.
%
%n appelle \emph{suite des meilleures approximations lin\'eaires} la suite d'\'el\'ements non nuls de $\mathbb{Z}^d$ de normes strictement croissante, not\'ee $(\Delta_n)_{n \in \mathbb{N}}$, d\'efinie au signe pr\`es, telle que quelque soit $ |\Delta | <  | \Delta_n |$, on ait

%$$ \| \langle \Delta_n,\theta \rangle \| < \| \langle \Delta,\theta \rangle \|. $$
%(notons que $ \| \langle \Delta,\theta \rangle \| = \| \langle \Delta',\theta \rangle \|$ implique $\Delta = \pm \Delta'$).

%On appelle \emph{suite des meilleures approximations simultan\'ees} la suite strictement croissante d'entiers strictement positifs, not\'ee $(q_n) _{n \in \mathbb{N}}$, telle que quelque soit $|q | < |q_n|$, on ait
%$$ \| q_n\theta\| < \| q\theta \|. $$
%
%
%
%
On utilisera \'egalement les notations suivantes : pour $h$ r\'eel non nul $$ \eps_s(h)=\min_{q \in \mathbb{Q}, |q| \le |h|}\|q\theta\|$$ et $$ \eps_l(h)=\min_{x \in \mathbb{Z}^d, |x|\le |h|}\| \langle x,\theta \rangle \|.$$

La suite des meilleures approximations simultan\'ees de $\theta$ v\'erifie l'in\'egalit\'e suivante
\begin{equation}
\label{inegmeilleurapprox}
\frac{1}{2}q_{n+1}^{-1} < (q_n+q_{n+1})^{-1} \le ||q_n\theta|| \le q_{n+1}^{-\frac{1}{d}}.
\end{equation}
Ce lemme est bien connu en dimension $1$ et provient des propri\'et\'es du d\'evelop\-pe\-ment en fractions continues. En dimension sup\'erieure,  l'in\'egalit\'e de droite se montre avec le principe de Dirichlet. Une d\'emonstration de l'in\'egalit\'e de gauche est donn\'ee dans \cite{Bernard}.

%1.3.1
\subsubsection{Notions habituelles d'approximation diophantienne en dimension $d$.}
Soit $\tau$ un r\'eel positif ou nul.

On rappelle que $\theta$ est diophantien de type $\tau$ pour l'approximation simultan\'ee si
$$\inf_{q \ne 0} q^{\frac{1+\tau}{d}} \|q \theta \| > 0,$$
c'est-\`a-dire avec nos notations si $\inf_{q \ne 0} q^{\frac{1+\tau}{d}} \eps_s(q) > 0.$
\medskip

On rappelle que $\theta$ est diophantien de type $\tau$ pour l'approximation lin\'eaire si
$$\inf_{|\Delta| \ne 0} |\Delta|^{d(1+\tau)} \| \langle \Delta,\theta \rangle \| > 0 ,$$
c'est-\`a-dire si $\inf_{h \ne 0} h^{d(1+\tau)} \eps_l(h) > 0 .$

On note $\Omega^d_s(\tau)$ (respectivement $\Omega^d_l(\tau)$) l'ensemble des vecteurs diophantiens de type $\tau$ pour l'approximation simultan\'ee (respectivement pour l'approximation lin\'eaire).

Les deux types d'approximations sont li\'es par le th\'eor\`eme de transfert de Khintchine (\cite{Kintchine1}, voir \cite{Cassels} et \cite{Laurent1} en dimension $2$).
\begin{thm*}
\label{thmtransfertkintchine}
[Khintchine]
Pour tout $\tau \ge 0$,
$$\Omega^d_s \left( \frac{\tau}{(d-1)\tau+d}\right) \subset \Omega^d_l(\tau) \subset \Omega^d_s(d\tau).$$
\end{thm*}

En particulier, on a $\Omega^d_l(0)=\Omega^d_s(0)$. On appelle vecteurs de type constant les vecteurs de type $0$. Cette notion concerne donc les m\^emes vecteurs pour les approximations lin\'eaires et simultan\'ees.

%1.3.2
\subsubsection{Une autre notion de type diophantien}

Nous nous int\'eressons \`a une notion diff\'erente introduite par Jarnik dans \cite{Jarnik3}.
%o\`u il prolonge les travaux qu'il avait effectu\'e sur le th\'eor\`eme de transfert de Khintchine permettant de montrer qu'il est optimal (\cite{Jarnik1}, \cite{Jarnik2}).
On pourra voir aussi un article de Khintchine (\cite{Kintchine3}). Laurent \cite{Laurent1} reprend cette notion en vue d'un r\'esultat qui pr\'ecise le th\'eor\`eme de transfert de Khintchine.

\begin{de}
\label{defomega}
Soit $\tau \ge 0$.

On note $\Theta_l^d(\tau)$ l'ensemble des vecteurs $\theta$ de $\mathbb{R}^d$ tels que
$$\limsup_{h \rightarrow +\infty} h^{d(1+\tau)} \eps_l(h) >0,$$
et $\Theta^d_s(\tau)$ l'ensemble des vecteurs $\theta$ de $\mathbb{R}^d$ tel que
$$\limsup_{q \rightarrow +\infty} q^{\frac{1+\tau}{d}} \eps_s(q) >0.$$

\end{de}

Autrement dit $\theta$ appartient \`a $\Theta_l^d(\tau)$ s'il existe une constante $C$ tel que pour une infinit\'e d'entiers $n $, pour tout $ \Delta$ de $\mathbb{Z}^d$, v\'erifiant $0 < | \Delta| \le n$, on ait $ \| \langle \Delta ,\theta \rangle \| \ge Cn^{-(1+\tau)d}.$ De m\^eme, $\theta$ appartient \`a $\Theta_s^d(\tau)$ s'il existe une constante $C$ telle que pour une infinit\'e d'entiers strictement positifs $n$, pour tout tout entier $q$ avec $ <q \le n$, on ait $\| q\theta\| \ge Cn^{-\frac{(1+\tau)}{d}}.$

On a bien s\^ur, $\Omega^d_s(\tau) \subset \Theta_s^d(\tau)$ et $\Omega^d_l(\tau) \subset \Theta_l^d(\tau)$. Pour mieux pointer la diff\'erence entre ces d\'efinitions et les notions habituelles d'approximation diophantienne, r\'e\'ecrivons-les \`a l'aide des meilleures approximations. Si on a $q_n \le q < q_{n+1}$, alors $\eps_s(q) = \|q_n\theta\|$. D'o\`u $\theta \in \Theta^d_s(\tau)$ si et seulement si
$$\limsup_{n \rightarrow +\infty} q_{n+1}^{\frac{1+\tau}{d}} \|q_n\theta\| >0, $$
tandis que $\theta \in \Omega^d_s(\tau)$ si et seulement si
$$\liminf_{n \rightarrow +\infty} q_{n}^{\frac{1+\tau}{d}} \|q_n\theta\| > 0.$$
La diff\'erence est la m\^eme pour les approximations lin\'eaires.
%1.3.3
\subsubsection{Le cas $\Theta^d_s(0)$.}
D'apr\`es l'in\'egalit\'e (\ref{inegmeilleurapprox}), tout vecteur de $\mathbb{R}^d$ \`a coordonn\'ees rationnellement ind\'e\-pendantes modulo un appartient \`a $\Theta_s^d(d-1)$. Le cas de la dimension 1 est donc particulier puisque tout $\theta$ irrationnel appartient \`a $\Theta_s^1(0)$.

En dimension $d>1$, on montrera que ce n'est plus vrai (voir par exemple la partie 4). Toutefois Chevallier a montr\'e dans \cite{Chevallier1} que pour presque tout $\theta$, \begin{large}                                                                                                                                                                               \end{large}on a $$ \limsup q_{n+1} || q_n \theta||^d >0.$$ C'est-\`a-dire que $\Theta^d_s(0)$ est de mesure $1$.

%1.4
\subsection{R\'esultats}
Ces notions diophantiennes moins usuelles vont nous permettre de donner une quasi-caract\'erisation des vecteurs $\theta$ pour lesquels la translation $T_\theta$ a ou n'a pas la propri\'et\'e du logarithme.
\begin{thm}
\label{thm1}
(i) \ Si $\theta$ appartient \`a $\Theta^d_s(0)$ alors la translation $T_{\theta}$ poss\`ede la propri\'et\'e du logarithme.
\\
(ii) \ S'il existe $\tau > 0$ tel que $\theta$ n'appartienne pas \`a $\Theta^d_s(\tau)$ alors $T_{\theta}$ ne poss\`ede pas la propri\'et\'e du logarithme.
\end{thm}
On a l'analogue pour l'approximation lin\'eaire.
\begin{thm}
\label{thm2}
(i) \ Si $\theta$ appartient \`a $\Theta^d_l(0)$ alors la translation $T_{\theta}$ poss\`ede la propri\'et\'e du logarithme.
\\
(ii) \ S'il existe $\tau > 0$ tel que $\theta$ n'appartienne pas \`a $\Theta^d_l(\tau)$ alors $T_{\theta}$ ne poss\`ede pas la propri\'et\'e du logarithme.
\end{thm}
%Remarquons toutefois qu'il existe des points du tore n'appartenant ni \`a $\Theta^d_s(0)$ ni \`a $\bigcap_{\tau>0} \Theta^d_s(\tau)$.
%mettre une note \`a la fin
Dans la partie 2, nous d\'emontrons le th\'eor\`eme \ref{thm1} \emph{(i)} et le th\'eor\`eme \ref{thm2} \emph{(ii)}. Pour finir de montrer les th\'eor\`emes \ref{thm1} et \ref{thm2} nous allons d\'emontrer dans la partie 3 une variante adapt\'ee \`a notre situation d'un th\'eor\`eme de transfert d\^u \`a Jarnik (\cite{Jarnik3}) :
%(semblable au th\'eor\`eme de transfert de Kintchine) entre les $\Theta_s$ et les $\Theta_l$. %
\begin{thm}
\label{thmtransfert}
Pour tout $\tau \ge 0$,
$$\Theta^d_s \left( \frac{\tau}{(d-1)\tau+d}\right) \subset \Theta^d_l(\tau) \subset \Theta^d_s(d\tau).$$
\end{thm}

En effet, cela implique $\Theta^d_s(0)=\Theta_s^d(0)$ et donc que les \'enonc\'es $(i)$ des th\'eor\`emes \ref{thm1} et \ref{thm2} sont \'equivalents. D'apr\`es l'inclusion de droite la condition $(ii)$ du th\'eor\`eme \ref{thm1} implique la condition $(ii)$ du th\'eor\`eme \ref{thm2}.

En dimension $1$, les parties $(i)$ des th\'eor\`emes \ref{thm1} et \ref{thm2} montrent en particulier que toute translation irrationnelle poss\`ede la propri\'et\'e du logarithme. En dimension sup\-\'erieure, d'apr\`es le r\'esultat de Chevallier cit\'e dans la partie 1.3.3, nous obtenons que pour presque tout $\theta $ la translation $T_\theta$ poss\`ede la propri\'et\'e du logarithme.

Dans la derni\`ere partie, pour $d=2$, nous construirons des vecteurs dont nous contr\^olons les approximations diophantiennes. Cela permet de montrer le th\'eor\`eme suivant :
\begin{thm}
\label{thm4}

(i) \ Il existe des vecteurs $\theta$ dans l'intersection de tous les $\Omega_s^2(\tau)$ pour $\tau$ strictement positif, pour lesquels $T_\theta$ ne poss\`ede pas la propri\'et\'e du logarithme.

(ii) \ Il existe des vecteurs $\theta$ pour lesquels $T_\theta$ poss\`ede la propri\'et\'e du logarithme et qui n'appartiennent \`a aucun $\Omega_s^2(\tau)$.
\end{thm}
\noindent
D'apr\`es le th\'eor\`eme de transfert de Khintchine (\ref{thmtransfertkintchine}), on a les m\^emes \'enonc\'es pour les $ \Omega_l^2(\tau)$.

\subsection{Propri\'et\'e du logarithme et loi du logarithme}

%Soit  $d_n(x,x_0) = \min_{1 \le i \le n} | T_\theta^i(x)-x_0|$.
Nous donnons le lien entre la propri\'et\'e du logarithme et la loi du logarithme et nous en d\'eduisons des r\'esultats. On suit la d\'efinition g\'en\'erale de \cite{Galatolo1} et on l'applique aux translations sur le tore de dimension $d$.
\begin{de} On dit que la translation $T_\theta$ v\'erifie la loi du logarithme si pour presque tout $x$,
$$\limsup_{n \to \infty}  \frac{-\log{\| T_\theta^n(x)\|}}{\log{n}} = \frac{1}{d}.$$
\end{de}

\begin{lem}
\label{prop-loi}
 La propri\'et\'e du logarithme implique la loi du logarithme.
\end{lem}

Remarquons que si la suite $( B( 0 , n^{-{\frac{1}{d}}} ) )_{n \in \mathbb{N}^*}$ est de Borel-Cantelli, alors pour presque tout $x$ il existe alors une infinit\'e de $n$ tels que $\| T_\theta^n(x)\| < n^{-{\frac{1}{d}}}$, d'o\`u
$$\limsup_{n \to \infty}  \frac{-\log{\| T_\theta^n(x)\|}}{\log{n}} \ge \frac{1}{d}.$$

Inversement, \'etant donn\'e un r\'eel $\delta$ strictement positif, si $$\limsup_{n \to \infty}  \frac{-\log{\| T_\theta^n(x)\|}}{\log{n}} > \frac{1}{\delta},$$ on a, de la m\^eme mani\`ere, que la suite $( B( 0 , n^{-{\frac{1}{\delta}}} ) )_{n \in \mathbb{N}^*}$ est de Borel-Cantelli. Alors la somme des mesures des boules est n\'ecessairement infinie, et $\delta$ ne peut \^etre strictement inf\'erieur \`a $d$. On a donc toujours $$ \limsup_{n \to \infty}  \frac{-\log{\| T_\theta^n(x)\|}}{\log{n}} \le \frac{1}{d}.$$

\qed
\medskip

L'article de Galatolo et Peterlongo (\cite{Galatolo1}) donne, dans le cas des translations sur le tore de dimension $2$, des contre-exemples \`a la loi du logarithme, qui sont donc \'egalement des contre-exemples \`a la propri\'et\'e du logarithme.
\medskip

\noindent
\textbf{Remarque}

Les th\'eor\`emes 1 et 2 restent vrais si on remplace ``propri\'et\'e du logarithme'' par ``loi du logarithme''. Pour les parties $(i)$, cela r\'esulte du lemme. Nous prouverons que c'est \'egalement vrai pour la partie $(ii)$ du th\'eor\`eme 2 (corollaire \ref{corloilog}) et donc aussi pour la partie $(ii)$ du th\'eor\`eme 1 par le th\'eor\`eme de transfert.

%2
\section{Crit\`eres sur $\theta$ pour que $T_{\theta}$ poss\`ede ou non la propri\'et\'e du logarithme}
%
%2.1
%
\subsection{D\'emonstration du th\'eor\`eme $1 (i)$ }

Soit $\theta \in \Theta_s^d(0)$ et $(q_n)_{n \in \mathbb{N}}$ la suite des meilleures approximations simultan\'ees de $\theta$, il existe donc $C>0$, tel que l'on ait pour une infinit\'e de $n$, $$q_{n+1}^\frac{1}{d} \| q_n \theta \| \ge C.$$

Soit $q \in \mathbb{N}$. On choisit un $n$ tel que $q_{n+1} \ge q$ et $q_{n+1}^\frac{1}{d} \| q_n \theta \| \ge C$ et on note
$$ U = \bigcup_{l=q}^{q+q_{n+1}-1} T_{\theta}^{-l} B(0,\frac{1}{l^{\frac{1}{d}}}).$$
Les rayons de ces $q_{n+1}$ boules sont sup\'erieurs \`a $(2q_{n+1})^{-\frac{1}{d}}$.
Les $q_{n+1}$ points $T^{-j}_\theta0$, avec $q \le j \le q +q_{n+1}-1$ sont \`a distance les uns des autres d'au moins $\| q_n \theta \|$ et $ \| q_n \theta \| \ge Cq_{n+1}^{-\frac{1}{d}}$. Donc $U$ contient $q_{n+1}$ boules disjointes de rayons $cq_{n+1}^{-\frac{1}{d}}$, avec $ c=\min \left( \frac{C}{2}, (\frac{1}{2})^{\frac{1}{d}} \right)$.

Il en r\'esulte
$$\mu(U) \ge q_{n+1} \frac{2^d {c}^d}{q_{n+1}} \ge 2^d{c}^d.$$
On a donc, quel que soit $q$ entier, $\mu \left( \cup_{l \ge q} T_{\theta}^{-l} B(0,\frac{1}{l^{\frac{1}{d}}}) \right) \ge 2^d{c}^d >0,$ d'o\`u
$$ \mu \left( \limsup T_{\theta}^{-l} B(0,\frac{1}{l^{\frac{1}{d}}}) \right) >0 .$$
D'apr\`es le lemme \ref{borelcantelli}, la suite de boules est de Borel-Cantelli.
 %\vspace*{-1cm}
%
%
%
%
%2.2
\subsection{Un crit\`ere pour que $T_{\theta}$ ne poss\`ede pas la propri\'et\'e du logarithme }

Nous donnons une condition suffisante pour ne pas poss\'eder la propri\'et\'e du logarithme. L'id\'ee de la d\'emonstration est que lorsque les approximations lin\'eaires sont suffisamment bonnes, les \'el\'ements de l'orbite vont \^etre assez proches d'un hyperplan pour que la mesure occup\'ee par l'union  des boules soit petite. Nous en d\'eduirons le th\'eor\`eme $2$ \emph{(ii)}

\begin{thm}
\label{thmprecis}
S'il existe une suite de vecteurs \`a coefficients entiers, de normes strictement croissantes, $(X_n)_{n \in \mathbb{N}}$ telle que
\begin{equation} \label{theoreme5} \sum |X_{n+1}|^{\frac{d}{d+1}}\|\langle X_n, \theta \rangle \|^{\frac{1}{d+1}} < \infty ,\end{equation}
alors la translation $T_{\theta}$ ne poss\`ede pas la propri\'et\'e du logarithme.
\end{thm}
Supposons la condition (\ref{theoreme5}) v\'erifi\'ee. Notons $\eps_{n}=\|\langle X_n, \theta \rangle \|$ et $B_n=B(0,n^{-\frac{1}{d}}) $.  Nous devons montrer que $\mu \left( \limsup T_\theta^{-n}B_n \right) = 0$. Pour cela, il suffit de montrer qu'il existe une suite de r\'eels strictements positifs $(L_n)_{n >0}$ tendant vers l'infini telle que
\begin{equation}
\label{conditionthm5}
 \sum_n \mu \left( \bigcup_{L_n \le l < L_{n+1}} T_{\theta}^{-l}B_l \right) < + \infty.
\end{equation}
Etant donn\'ee une suite $(L_n)$ de r\'eels strictements positifs, soit $U_n=\bigcup_{L_n \le l < L_{n+1}} T_{\theta}^{-l}B_l$. Soient $x$ un point de $U_n$ et $k$ un entier avec $L_n \le k < L_{n+1}$, tel que $ x \in T_{\theta}^{-k}B_k$, c'est-\`a-dire  $\| x + k\theta \| \le L_n^{-\frac{1}{d}}$. On a alors
\begin{eqnarray}
\| \langle X_n,x \rangle \| &\le& k \| \langle X_n,\theta \rangle \| + \| \langle X_n,x + k\theta \rangle \|, \nonumber \\
&\le& L_{n+1}\eps_n+d \| x + k\theta \| \, |X_n|, \nonumber \\
 &\le& L_{n+1}\eps_n + d L_n^{-\frac{1}{d}} |X_n| \nonumber .
\end{eqnarray}

L'application $\tilde{X}_n$ de $\mathbb{T}^d$ dans $\mathbb{T}$ qui envoie $x$ sur $ \langle X_n,x \rangle \mod 1$ est un morphisme surjectif de groupes compacts. L'image de $\mu$ par $\tilde{X}_n$ est donc la mesure de Lebesgue sur $\mathbb{T}$ et il r\'esulte de l'in\'egalit\'e pr\'ec\'edente que $$\mu(U_n) \le 2 \left( L_{n+1}\eps_n + d L_n^{-\frac{1}{d}} |X_n| \right).$$

\noindent
Il nous suffit donc de construire une suite $(L_n)_{n \in \mathbb{N}}$ tendant vers l'infini telle que
\begin{equation}\label{deux} \sum_{n} \left( L_{n+1} \eps_{n} + |X_n| L_n^{-\frac{1}{d}} \right) < \infty.\end{equation}
On pose $L_{n} = |X_{n}|_{ }^{\frac{d}{d+1}} \eps_{n-1}^{-\frac{d}{d+1}}$. Cette suite tend bien vers l'infini et
$$L_{n} \eps_{n-1} = |X_{n}| L_{n}^{-\frac{1}{d}} = |X_{n}|^{\frac{d}{d+1}}\eps_{n-1}^{\frac{1}{d+1}}.$$
D'o\`u l'hypoth\`ese (\ref{theoreme5}) entra\^ine (\ref{deux}), ce qui conclut la d\'emonstration du th\'eor\`eme \ref{thmprecis}.
%Notre condition suffisante se r\'e\'ecrit donc
%
%$$\sum_{n \ge 1} |X_{n}|^{\frac{d}{d+1}}\eps_{n-1}^{\frac{1}{d+1}} < \infty .$$
%
%2.3

\subsection{D\'emonstration du th\'eor\`eme $2$ \emph{(ii)}}
Nous donnons des conditions \'equivalentes \`a la condition (\ref{theoreme5}) du th\'eor\`eme \ref{thmprecis}.
\begin{lem}\label{equivalence1}
Les propri\'et\'es suivantes sont \'equivalentes :

(i) Il existe $(X_{n})_{n \in \mathbb{N}}$, une suite de vecteurs non nuls \`a coefficients entiers de normes strictement croissantes, telle que
$\sum |X_{n+1}|^{\frac{d}{d+1}} \| \langle X_n,\theta \rangle \|^{\frac{1}{d+1}} < \infty.$

(ii) Il existe une sous-suite $(\Delta_{\phi(n)})_{n \in \mathbb{N}}$ de la suite des meilleures approximations lin\'eaires telle que $\sum |\Delta_{\phi(n+1)}|^{\frac{d}{d+1}} \| \langle \Delta_{\phi(n)},\theta \rangle \|^{\frac{1}{d+1}} < \infty.$

(iii) $\sum (2^{nd} \eps_l(2^n))^{\frac{1}{d+1}} < \infty.$

(iv) $\sum_{k \in \mathbb{N}^*} k^{-1}(k^d \eps_l(k))^{\frac{1}{d+1}} < \infty.$
%(v) On a $\int_1^\infty (\frac{ \eps_l(h)}{h})^{\frac{1}{d+1}} dh< \infty.$.
\end{lem}
Montrons d'abord qu'on d\'eduit le th\'eor\`eme \ref{thm2} \emph{(ii)} du lemme \ref{equivalence1}. Supposons qu'il existe $\tau$ strictement positif tel que $\theta$ n'appartienne pas \`a $\Theta_l^d(\tau)$. Alors pour tout $k$ assez grand, $\eps_l(k) < k^{-d(1 + \tau)}$. On a donc
$$k^{-1}(k^d \eps_l(k))^{\frac{1}{d+1}} \le k^{-1-\frac{d\tau}{d+1}}.$$
La propri\'et\'e \emph{(iv)} est v\'erifi\'ee puisque $\tau$ est strictement positif. D'apr\`es le lemme, l'hypoth\`ese du th\'eor\`eme \ref{thmprecis} est \'egalement v\'erifi\'ee et donc $T_{\theta}$ ne poss\`ede pas la propri\'et\'e du logarithme.
\par\medskip\noindent
\textbf{D\'emonstration du lemme \ref{equivalence1}}
%et $\eps_n = \| \langle \Delta_n,\theta \rangle \| $

Nous montrons tout d'abord que $(i)$ implique $(ii)$. Soit $(X_n)_{n \in \mathbb{N}}$ une suite de vecteurs v\'erifiant la condition $(i)$. On consid\`ere la suite $(\Delta_n)_{n \in \mathbb{N}}$ des meilleures approximations lin\'eaires. On d\'efinit la suite d'entiers $(\phi(n))_{n \in \mathbb{N}}$ par la relation $$|\Delta_{\phi(n)}| \le |X_{n}| < | \Delta_{\phi(n)+1}|.$$
D'apr\`es la d\'efinition des meilleures approximations lin\'eaires, on a pour tout $n$,
$\| \langle \Delta_{\phi(n)},\theta \rangle \| \le \| \langle X_n,\theta \rangle \|$. 
D'o\`u $$\sum \left(|\Delta_{\phi(n+1)}|^d \| \langle \Delta_{\phi(n)},\theta \rangle \|\right)^{\frac{1}{d+1}} < \sum \left( |X_{n+1} |^d \| \langle X_n,\theta \rangle \| \right)^{\frac{1}{d+1}} < + \infty .$$

La suite $(\phi(n))_{n \in \mathbb{N}}$ n'est pas n\'ecessairement strictement croissante, mais on se ram\`ene sans difficult\'e \`a ce cas.
\medskip

Nous montrons maintenant que $(ii)$ implique $(iii)$. Soit $(\Delta_{\phi(n)})$ une sous-suite des meilleures approximations lin\'eaires v\'erifiant la condition $(ii)$.
%$$\sum |\Delta_{\phi(n+1)}|^{\frac{d}{d+1}} \| \langle \Delta_{\phi(n)},\theta \rangle \|^{\frac{1}{d+1}} < \infty.$$
%
Pour $k \ge 0 $, notons  $h_k = |\Delta_{\phi(k)}| $. Pour tout entier $n$ tel que $h_k \le 2^n < h_{k+1},$ on a $\eps_l(2^n) \le \| \langle \Delta_{\phi(k)}, \theta \rangle \|,$ d'o\`u

\begin{eqnarray}
\sum_{h_k \le 2^n < h_{k+1}}(2^{nd}\eps_l(2^n))^{\frac{1}{d+1}} &\le& \| \langle \Delta_{\phi(k)}, \theta \rangle \|^{\frac{1}{d+1}} \sum_{h_k \le 2^n < h_{k+1}}(2^{nd})^{\frac{1}{d+1}} \nonumber \\
&\le& \| \langle \Delta_{\phi(k)}, \theta \rangle \|^{\frac{1}{d+1}} \frac{ (2 h_{k+1})^{\frac{d}{d+1}}}{2^{\frac{d}{d+1}}-1}. \nonumber
\end{eqnarray}
Donc,
$$\sum_{n \ge n_0} \left(2^{nd} \eps_l(2^n)\right)^{\frac{1}{d+1}} \le c \sum_{k \ge 0} h_{k+1}^{\frac{d}{d+1}} \| \langle \Delta_{\phi(k)}, \theta \rangle \|^{\frac{1}{d+1}},$$
o\`u  $n_0$ est un entier v\'erifiant $2^{n_0} \ge h_0$ et $c$ est une constante strictement positive.
\medskip

Maintenant, pour d\'eduire la propri\'et\'e $(i)$ de la propri\'et\'e $(iii)$, il suffit de remarquer qu'\`a $\eps_l(2^n)$ correspond un suite de droites $(X_{n})$, avec $|X_{n}| \le 2^n$, v\'erifiant
$$\| \langle X_n,\theta \rangle \| = \min_{|X| \le 2^n} \| \langle X,\theta \rangle \| = \eps_l(2^n) .$$
On a alors
$$|X_{n+1}|^{\frac{d}{d+1}} \| \langle X_n,\theta \rangle \|^{\frac{1}{d+1}} \le 2^{(n+1)\frac{d}{d+1}}\eps_l(2^n)^{\frac{1}{d+1}} = 2^{\frac{d}{d+1}} \left (2^{nd}\eps_l(2^n) \right )^{\frac{1}{d+1}}.$$

On se ram\`ene ensuite \`a une suite dont les normes sont strictement croissantes.
\medskip

Nous finissons en montrant que les propri\'et\'es $(iii)$ et $(iv)$ sont \'equivalentes. Remarquons que pour $n$ donn\'e si on a $2^n \le k < 2^{n+1}$, d'apr\`es les propri\'et\'es de meilleures approximations, $\eps_l(2^{n+1}) \le \eps_l(k) \le \eps_l(2^{n})$ et donc
$$\frac{1}{2^{n+1}}( 2^{nd}\eps_l(2^{n+1}))^{\frac{1}{d+1}} \le \frac{1}{k} (k^d \eps_l(k))^{\frac{1}{d+1}} \le \frac{1}{2^{n}}( 2^{(n+1)d}\eps_l(2^{n}))^{\frac{1}{d+1}}.$$
En sommant ces in\'egalit\'es pour $k$ compris entre $2^n$ et $2^{n+1}$
$$\frac{1}{2}( 2^{nd}\eps_l(2^{n+1}))^{\frac{1}{d+1}} \le \sum_{2^n \le k < 2^{n+1}} \frac{1}{k}(k^d \eps_l(k))^{\frac{1}{d+1}} \le ( 2^{(n+1)d}\eps_l(2^{n}))^{\frac{1}{d+1}}.$$
Puis on somme sur les entiers $n \ge 0 $,
$$2^{-(1+\frac{d}{d+1})}\sum_{n \ge 1} (2^{nd} \eps_l(2^n))^{\frac{1}{d+1}} \le \sum_{k \ge 1}\frac{1}{k} (k^d \eps_l(k))^{\frac{1}{d+1}} \le 2^{\frac{d}{d+1}} \sum_{n \ge 0} (2^{nd} \eps_l(2^n))^{\frac{1}{d+1}} .$$

%%
%

%2.4
\subsection{Un premier contre-exemple}
Donnons maintenant un premier contre-exemple \`a la propri\'et\'e du logarithme. Cet exemple s'inspire de l'id\'ee d'alterner les meilleurs approximations des coordonn\'ees de l'angle $\theta$. Cette id\'ee a \'et\'e utilis\'ee par Yoccoz pour d\'emontrer que la propri\'et\'e de Denjoy-Koksma n'\'etait plus vraie en dimension sup\'erieure \`a 1 (\cite{Yoccoz}). On la retrouve par exemple dans \cite{Fayad2} pour montrer que l'on peut construire des flots sp\'eciaux au-dessus de rotations sur le tore de dimension 2 qui ne soient pas m\'elangeants, dans la construction de contre-exemples \`a la loi du logarithme par Galatolo et Peterlongo (\cite{Galatolo1}) ou dans le contre-exemples que donne Chevallier (\cite{Chevallier2}) d'un point dont la trajectoire est ``mal r\'epartie". 
\medskip

Pour $d>1$, soit $\theta = (\theta_1,...,\theta_d)$. Pour $n \in \mathbb{N}$ et pour $1 \le i \le d$, on note $(q_{i,n})$ la suite des d\'enominateurs de la fraction continue de $\theta_i$ et on pose $X_{dn+i}=(0,...,q_{i,n},..,0)$. On a alors en particulier $$\| \langle X_{dn +i}, \theta \rangle \| \le \frac{1}{q_{i,n+1}}.$$
D'apr\`es le th\'eor\`eme \ref{thmprecis}, la translation $T_\theta$ n'a pas la propri\'et\'e du logarithme si ces suites v\'erifient
$$\sum_{n} \left( \left(\frac{q_{1,n}^d}{q_{d,n}} \right)^{\frac{1}{d+1}}+\left(\frac{q_{2,n}^d}{q_{1,n+1}} \right)^{\frac{1}{d+1}}+ ...+\left(\frac{q_{d,n}^d}{q_{d-1,n+1}} \right)^{\frac{1}{d+1}} \right) < \infty.$$
Cette condition est r\'ealis\'ee s'il existe une constante $c>1$ pour $n$ assez grand, tel que $$\left(\frac{q_{1,n}^d}{q_{d,n}} \right)^{\frac{1}{d+1}} \le n^{-c}$$
et pour $2 <i \le d,$ $$ \left(\frac{q_{i,n}^d}{q_{i-1,n+1}} \right)^{\frac{1}{d+1}} \le n^{-c}.$$
On en d\'eduit la proposition suivante
\begin{prop}
\label{propcontreexemple} Soit $d>1$. Soit $\theta = (\theta_1,...,\theta_d)$ et soit pour $1 \le i \le d$,  la suite $(q_{i,n})$ des d\'enominateurs des fractions continues de $\theta_i$. S'il existe $\delta > d+1$ tel que pour tout $n$ assez grand, $$q_{d,n} \ge q_{1,n}^d n^\delta \; \mbox{et} \; q_{i-1,n+1} \ge q_{i,n}^d n^\delta  \; \mbox{pour} \; 2 < i \le d,$$ alors $T_{\theta}$ ne poss\`ede pas la propri\'et\'e du logarithme.
\end{prop}
%c'est-\`a-dire sa suite des meilleures approximations

Il est ais\'e de construire des vecteurs $\theta_i$ v\'erifiant les relations de la proposition \ref{propcontreexemple}, en construisant par r\'ecurrence, simultan\'ement, les d\'eveloppements en fractions continues des $\theta_i$.

En dimension 2, on montre (voir \cite{Fayad2} chapitre 7) que quel que soit $\eps$, il existe des vecteurs $\theta \in \Omega^2_s(1+\eps)$ v\'erifiant ces conditions. Par contre, on peut montrer (voir \cite{Bernard}) que de tels $\theta$ ne peuvent avoir des types diophantiens plus petits.
%
%2.5
\subsection{Extensions du crit\`ere \`a des boules de rayon $n^{-\frac{1}{\delta}}$}
Nous \'elargissons dans cette section le cadre du probl\`eme. Remarquons que dans la d\'emonstration du th\'eor\`eme \ref{thmprecis} le fait que les rayons des boules soient \'egaux \`a $n^{-\frac{1}{d}}$ ne joue pas un r\^ole important. Soit$(X_n)$ une suite  de vecteurs \`a coefficients entiers. On consid\`ere les boules $B(0,r_n)$, avec $r_n $ d\'ecroissant et une suite $(L_n)$ de r\'eels strictement positifs tendant vers l'infini. Comme dans la d\'emonstration du th\'eor\`eme \ref{thmprecis}, on note  $\eps_{n}=\|\langle X_n, \theta \rangle \|$ et on d\'efinit $U_n=\bigcup_{L_n \le  l < L_{n+1}} T_{\theta}^{-l}B(0,r_l)$. On obtient de la m\^eme mani\`ere que si un point du tore $x$ appartient \`a $U_n$, alors
$$ \| \langle X_n,x \rangle \| \le L_{n+1}\eps_n + d r_{L_n} |X_n| .$$
D'o\`u
$$\mu(U_n) \le 2 ( L_{n+1} \eps_n+ dr_{L_n} |X_n|) .$$

En particulier, lorsque $r_n = n^{-\frac{1}{\delta}}$, en posant, $$L_{n} = |X_{n}|_{ }^{\frac{\delta}{\delta+1}} \eps_{n-1}^{-\frac{\delta}{\delta+1}},$$ et comme pour le th\'eor\`eme \ref{thmprecis}, on obtient

%de la m\^eme mani\`ere que $\mu(U_n) \le 2 |X_{n}|_{ }^{\frac{\delta}{\delta+1}} \eps_{n-1}^{\frac{1}{\delta+1}}$. On a donc
\begin{prop}
\label{prop-delta}
Soit $\delta > 0$. S'il existe une suite de vecteurs \`a coefficients entiers, de normes croissantes, $(X_n)_{n \in \mathbb{N}}$ telle que
$$\sum |X_{n+1}|^{\frac{\delta}{\delta+1}}\eps_{n}^{\frac{1}{\delta+1}} < \infty ,$$
o\`u $\eps_{n}=\|\langle X_n, \theta \rangle \|$, alors la suite des boules $B(0,n^{-\frac{1}{\delta}})$ n'est pas de Borel-Cantelli pour $T_\theta$.
\end{prop}
Remarquons que les suite de boules $B(0,n^{-\frac{1}{\delta}})$ avec $\delta < d$ ne peuvent pas \^etre de Borel-Cantelli, car la somme de leurs mesures est finie.

Dans la d\'emonstration du lemme \ref{equivalence1}, on peut sans difficult\'es remplacer $d$ par $\delta$, et on obtient de m\^eme
\begin{lem}
\label{lm-delta}
Soit $\delta>0$. Les propri\'et\'es suivantes sont \'equivalentes :

(i) Il existe $(X_{n})_{n \in \mathbb{N}}$, une suite de vecteurs \`a coefficients entiers de normes strictement croissantes, telle que
$\sum |X_{n+1}|^{\frac{\delta}{\delta+1}} \| \langle X_n,\theta \rangle \|^{\frac{1}{\delta+1}} < \infty.$

(ii) $\sum_{k \in \mathbb{N}^*} k^{-1}(k^\delta \eps_l(k))^{\frac{1}{\delta+1}} < \infty.$
%(v) On a $\int_1^\infty (\frac{ \eps_l(h)}{h})^{\frac{1}{d+1}} dh< \infty.$.
\end{lem}
Il en r\'esulte
\begin{prop}
\label{thmdelta}
Soit $\delta \ge d$. Si $\theta$ n'appartient pas \`a $\Theta_l^d(\tau)$ pour un $\tau>\frac{\delta}{d}-1$, les suite de boules $B(0,n^{-\frac{1}{\delta}})$ ne sont pas de Borel-Cantelli pour $T_\theta$.
\end{prop}

En effet, supposons qu'il existe $\tau$ avec $\tau>\frac{\delta}{d}-1$ tel que $\theta$ n'appartienne pas \`a $\Theta_l^d(\tau)$. Alors pour tout $k$ assez grand, $\eps_l(k) < k^{-d(1 + \tau)}$. On a donc
$$k^{-1}(k^\delta \eps_l(k))^{\frac{1}{\delta+1}} \le k^{-1-\frac{ d+d\tau-\delta}{d+1}}$$
On a $d+d\tau-\delta > 0$, donc la propri\'et\'e \emph{(ii)} du lemme \ref{lm-delta} est v\'erifi\'ee et la proposition \ref{prop-delta} s'applique.

\begin{cor}
\label{corloilog}
S'il existe $\tau > 0$, tel que $\theta$ n'appartienne pas \`a $\Theta^d_l(\tau)$, alors la translation $T_\theta$ ne v\'erifie pas la loi du logarithme.
\end{cor}

En choisissant $ d < \delta < d( 1 + \tau)$, la suite de boules $B(0,n^{-\frac{1}{\delta}})$ n'est pas de Borel-Cantelli d'apr\`es la proposition \ref{thmdelta}. D'apr\`es la d\'emonstration du lemme \ref{prop-loi}, on ne peut pas alors avoir
$\limsup_{n \to \infty}  \frac{-\log{\| T_\theta^n(x)\|}}{\log{n}} > \frac{1}{\delta},$ d'o\`u
$$  \limsup_{n \to \infty}\frac{-\log{\| T_\theta^i(x)\|}}{\log{n}} < \frac{1}{d}.$$
%\limsup_{n \to \infty}

%
%En particulier il existe $\delta>d$ tel que $\tau>\frac{\delta}{d}-1$, d'o\`u d'apr\`es la proposition $2.5$,
%
%3
\section{Relation de transfert entre les $\Theta_l^d(\tau)$ et les $\Theta_s^d(\tau)$ }
\subsection{D\'emonstration du th\'eor\`eme de transfert}

Nous montrons le th\'eor\`eme \ref{thmtransfert} en en donnant une version plus pr\'ecise. Cela compl\'etera la d\'emonstration des th\'eor\`emes \ref{thm1} et \ref{thm2}. Nous en d\'eduirons \'egalement une variante du  th\'eor\`eme \ref{thmprecis} avec une condition portant sur les approximations simultan\'ees.

\begin{thm}
\label{thmtransfertprecis} (i) Quels que soient $h>0$  et $d$ entier non nul on a
$$\eps_l(h) \le \frac{1}{Ch^{d-1}} \, \eps_s(Ch^d),$$
o\`u $C=\frac{1}{2(d+1)}.$

(ii) Soient $\tau>0$ et $\eta>0$. Si $\limsup h^{d(1+\tau)} \eps_l(h) < \eta $, alors
$$\limsup q^{\frac{1+\tau}{d +(d-1)\tau}}\eps_s(q) \le C'\eta^{\frac{1}{d^2+d(d-1)\tau}},$$
o\`u $C'$ est une constante strictement positive.

\end{thm}

Le th\'eor\`eme \ref{thmtransfert} en r\'esulte bien. Soit en effet $\tau \ge 0$, en posant $q=Ch^d$, l'in\'egalit\'e $(i)$ s'\'ecrit $$h^{d(1+\tau)}\eps_l(h) \le c_1q^{\frac{1+d\tau}{d}}\eps_s(q),$$
o\`u $c_1$ est une constante strictement positive. Ce qui implique que $\Theta_l^d(\tau) \subset \Theta_s^d(d\tau).$

D'autre part, si $\theta$ n'appartient pas \`a $\Theta^d_l(\tau)$ alors $\limsup h^{d(1+\tau)} \eps_l(h)=0$ et d'apr\`es $(ii)$, $$\limsup q^{\frac{1+\tau}{d +(d-1)\tau}}\eps_s(q)=0.$$ Or on a $\frac{1+\tau}{d +(d-1)\tau} = \frac{1}{d}\left ( 1+ \frac{\tau}{d +(d-1)\tau} \right ) $, donc $\theta$ n'appartient pas \`a $\Theta^d_s \left( \frac{\tau}{(d-1)\tau+d}\right)$.
\medskip

%Si $\theta$ n'appartient pas \`a l'intersection des $\Theta_s^d(\tau)$ alors $\theta$ n'appartient pas \`a l'intersection des $\Theta_l^d(\tau)$
%Si $\theta$ appartient \`a $\Theta_l^d(0)$ alors $\theta$ appartient \`a $\Theta_s^d(0)$.
%(i) $\bigcap_{\tau >0} \Theta_l^d(\tau) \subset \bigcap_{\tau >0} \Theta_s^d(\tau)$ .
%

La preuve du th\'eor\`eme reprend la d\'emonstration par Khintchine (\cite{Kintchine3}) du th\'eor\`eme de transfert de Jarnik. A partir d'un lemme de Minkoswki, Khintchine montre le lemme ci-dessous. On pourra \'egalement voir un article d'Apfelbeck (\cite{Apf}) donnant une am\'elioration du th\'eor\`eme de transfert de Jarnik et qui reprend ce lemme clef de Khintchine.

%

%\emph{Si $f_1,...,f_n$ sont des formes lin\'eaires \`a coefficients r\'eels et $c_1,...,c_n$ de r\'eels strictement positifs tels que $|\det{(f_1,...,f_n)}| = \prod_{i=1}^n c_i$, il existe $x \in \mathbb{Z}^n \setminus \{0\}$ v\'erifiant $|f_1(x)| \le c_1$ et $|f_i(x)| < c_i $ pour $2\le i \le d$}.

%
\begin{lem}
\label{lemmemagique}
[Khintchine]
Soient $f_1, ... , f_n$ et $g_1, ... , g_n$ des formes lin\'eaires sur $\mathbb{R}^n$, avec $\det(g_1,...,g_n) = \lambda$, telles que la forme bilin\'eaire sur $\mathbb{R}^{n }\times \mathbb{R}^{n}$, $(u,v) \to \sum_{i=1}^{n} f_i(u)g_i(v)$ soit \`a coefficients entiers, soient $t_1,..., t_n$ des r\'eels strictement positifs.
S'il existe $u \in \mathbb{Z}^n \setminus\{0\}$ tel que $|f_i(u)| \le t_i$, pour $1 \le i \le n$ et $f_i(u) \ne 0$ pour au moins un $i$, alors il existe $v \in \mathbb{Z}^n \setminus\{0\}$ tel que pour $1 \le k \le n$,
$$| g_k(v)| \le (2n \lambda )^{\frac{1}{n-1}}\frac{(\prod_{i=1}^n t_i)^{\frac{1}{n-1}}}{t_k} .$$
\end{lem}
\noindent
\textbf{D\'emonstration du th\'eor\`eme \ref{thmtransfertprecis}}

On applique le lemme, avec $n=d+1$, aux formes lin\'eaires sur $\mathbb{R}^{d+1}$ d\'efinies par
\begin{equation*}
\begin{cases}
f_i(u) = u_i - \theta_i u_{d+1}, \; ${\footnotesize $1 \le i < d$ }$\\
f_{d+1}(u) = u_{d+1} \\
g_i(v) = v_i, \; ${\footnotesize $1 \le i \le d$ }$\\
g_{d+1}(v) = \sum_{i=1}^d v_i \theta_i + v_{d+1}

\end{cases}
\end{equation*}

On trouve $\sum_{i=1}^{d+1} f_i(u)g_i(v) = \sum_{i=1}^{d+1} u_iv_i$, il s'agit donc bien d'une forme bilin\'eaire \`a coefficients entiers. On a $\det(g_1,...,g_{d+1}) = 1$ et $\det(f_1,...,f_{d+1}) = 1$.
\medskip

On commence par montrer la partie  $(i)$. Soient $q$ un r\'eel strictement positif et $\eps = \eps_s(q)$. On choisit $u = (u_1,...,u_{d+1})$ dans $\mathbb{Z}^{d+1}$ tel que $|u_{d+1}| < q $ et $\|u_{d+1} \theta \| = \max_{1 \le i \le d} |u_{d+1}\theta_i-u_i | = \eps$. Alors $|f_i(u)| \le \eps$, pour $1 \le i \le d$, et $|f_{d+1}(u)| \le q$.

D'apr\`es le lemme, il existe $v=(v_1,...,v_{d+1})$ non nul appartenant \`a $\mathbb{Z}^{d+1}$ tel que, en notant $\Delta=(v_1,...,v_d)$ on ait,
\begin{equation*}
\begin{cases}
|\Delta| = \max_{1 \le i \le d} |v_i| = \max_{1 \le i \le d} |g_i(v)| \le c \,  \displaystyle{\frac{(\eps^dq)^{\frac{1}{d}}}{\eps}} = cq^{\frac{1}{d}} \\
\| \langle \Delta, \theta \rangle \| = | \sum_{i=1}^d v_i \theta_i + v_{d+1} | = |g_{d+1}(v)| \le c \, \displaystyle{\frac{(\eps^dq)^{\frac{1}{d}}}{q}} = c q^{-{\frac{d-1}{d}}} \eps
\end{cases}
\end{equation*}
o\`u $c=\left(2(d+1)\right)^{\frac{1}{d}}$.

Il en r\'esulte que
$$\eps_l(cq^{\frac{1}{d}}) \le c  q^{-(1- \frac{1}{d})} \eps.$$

Pour $h> 0$ donn\'e, en posant $q=Ch^d$, c'est-\`a-dire $h=cq^{\frac{1}{d}}$, on trouve bien
$$ \eps_l(h) \le \frac{1}{Ch^{d-1}} \ \eps_s(Ch^d).$$

%Pour montrer le th\'eor\`eme $6 (ii)$, nous allons passer, en reprenant une id\'ee d'Apfelbeck, par le lemme suivant :

%\begin{lem}
%S'il existe deux suite de r\'eels strictement positifs $(K_m)$ et $(H_m)$ tels que pour tout m entier, pour tout $h > H_m$ on ait $h^{d(1+\tau)} \eps_l(h) < K_m$, alors il existe une suite $(Q_m)$ tels que pour tout $q > Q_m$, on ait
%$$q^{\frac{1+\tau}{d +(d-1)\tau}}\eps_s(q) \le c'K_m^{\frac{1+(d-2)\tau}{d^2+d(d-1)\tau}},$$
%avec $c'>0$ une constante d\'ependant uniquement de $d$ et $\tau$.
%\end{lem}

%Ce lemme implique, en faisant tendre la suite $(K_m)$ vers $0$, que si $\theta$ n'appartient pas \`a $\Theta^d_l(\tau)$ alors $\theta$ n'appartient pas \`a $\Theta^d_s \left( \frac{\tau}{(d-1)\tau+d}\right)$, ce qui conclut la d\'emonstration du th\'eor\`eme $6 (ii)$.

On montre maintenant la partie $(ii)$ du th\'eor\`eme. On utilise le lemme \ref{lemmemagique} avec les m\^emes formes lin\'eaires que dans la d\'emonstration de la partie $(i)$, mais on inverse les r\^oles des $f_i$ et des $g_i$.

Etant donn\'e $\tau>0$ et $\eta>0$, soit $h>0$ tel que $h^{d(1+\tau)} \eps_l(h) < \eta$, et soit $v$ un vecteur non nul de $\mathbb{Z}^d$ tel que $|g_i(v)| = |v_i| \le |h|$, pour $1\le i \le d$, et $$|g_{d+1}(v)|  = \eps_l(h) \le \frac{\eta}{h^{d(1 + \tau)}} \; \cdot$$

D'apr\`es le lemme \ref{lemmemagique}, il existe un vecteur $u$ non nul de $\mathbb{Z}^{d+1}$ tel que pour $1 \le i \le d$, on a $$|f_i(u)| = | u_{d+1}\theta - u_i | \le  \frac{c}{h}\left ( \frac{h^d \eta}{h^{d(1 + \tau)}} \right )^{\frac{1}{d}} = c \ \frac{\eta^{\frac{1}{d}}}{h^{(1 + \tau)}} \; ,$$
et
$$|f_{d+1}(u)| = |u_{d+1}| \le c \ \frac{h^{d(1 + \tau)}}{\eta}\left ( \frac{h^d \eta}{h^{d(1 + \tau)}} \right )^{\frac{1}{d}} = c \eta^{\frac{1}{d}-1}h^{d + (d-1)\tau}.$$

En posant $q=c \eta^{\frac{1}{d}-1}h^{d + (d-1)\tau}$, on obtient $$\eps_s(q) \le c \eta^{\frac{1}{d}}h^{-(1 + \tau)}.$$

\noindent
Comme $$q^{\frac{1+\tau}{d +(d-1)\tau}} = \left(c\eta^{\frac{1}{d}-1}\right)^{\frac{1+\tau}{d+(d-1)\tau}} h^{(1+\tau)},$$ 

$$q^{\frac{1+\tau}{d +(d-1)\tau}}\eps_s(q) \le\left(c\eta^{\frac{1}{d}-1}\right)^{\frac{1+\tau}{d+(d-1)\tau}} c \eta^{\frac{1}{d}} =  C'\eta^{\frac{1}{d^2+d(d-1)\tau}},$$

\noindent
o\`u $C'$ est une constante strictement positive.

Sous l'hypoth\`ese de $(ii)$, on peut choisir $h$ et donc $q$ arbitrairement grands v\'erifiant ces in\'egalit\'es, ce qui conclut la d\'emonstration du th\'eor\`eme \ref{thmtransfertprecis}.

%
%3.2
\subsection{Crit\`ere en approximations simultan\'ees pour ne pas avoir la propri\'et\'e du logarithme}

\begin{prop}
\label{prop32}
S'il existe une suite d'entiers strictement croissante $(q_n)_{n \in \mathbb{N}}$, telle que $$\sum_{n \ge 0} (q_n^{\frac{1}{d}} \|q_{n-1} \theta \|)^{\frac{1}{d+1}} < \infty,$$ alors $T_\theta$ ne poss\`ede pas la propri\'et\'e du logarithme.
\end{prop}

Commen\c cons par r\'e\'ecrire la condition \emph{(iv)} du lemme \ref{equivalence1} sous la forme
$$ \int_1^\infty \frac{1}{h}(h^d \eps_l(h))^{\frac{1}{d+1}} dh < \infty.$$
D'apr\`es le th\'eor\`eme \ref{thmtransfertprecis}, on a
$$ \int_1^\infty \frac{1}{h} (h^d \eps_l(h))^{\frac{1}{d+1}} dh \le \frac{1}{C^{\frac{1}{d+1}}} \int_1^\infty \frac{1}{h} (h \eps_s(Ch^d))^{\frac{1}{d+1}} dh.$$
On applique le changement de variable $q=Ch^{d}$,
$$ \int_1^\infty \frac{1}{h} (h^d \eps_l(h))^{\frac{1}{d+1}} dh  \le \frac{1}{d}\frac{1}{C^{\frac{1}{d}}} \int_{C}^\infty \frac{1}{q}(q^{\frac{1}{d}} \eps_s(q))^{\frac{1}{d+1}}dq,$$
%\le c_2 \int_{C}^\infty \frac{1}{q^{\frac{1}{d}}}(q^{\frac{1}{d}} \eps_s(q))^{\frac{1}{d+1}}\frac{1}{q}q^{\frac{1}{d}}dq

Donc si $\sum_{q \in \mathbb{N}^*} q^{-1}\left (q^{\frac{1}{d}} \eps_s(q) \right )^{\frac{1}{d+1}}< \infty $, alors $T_\theta$ ne poss\`ede pas la propri\'et\'e du logarithme.

On montre que cette condition est \'equivalente \`a la condition de l'\'enonc\'e de la m\^eme mani\`ere que dans la d\'emonstration du lemme \ref{equivalence1}.
%
%%4
\section{Construction, en dimension $2$, d'exemples et de contre-exemples \`a la propri\'et\'e du logarithme}
Dans $\cite{Laurent1}$, Laurent d\'emontre l'optimalit\'e d'une in\'egalit\'e de transfert en construisant un vecteur $\theta$ dont il contr\^ole les approximations diophantiennes. Nous allons reprendre sa construction en l'adaptant \`a nos d\'efinitions.
\begin{thm}
\label{thm7}
Etant donn\'e $(a_n)_{n \in \mathbb{N}}$ et $(h_n^\circ)_{n \in \mathbb{N}}$ deux suites d'entiers strictement positifs, avec $a_n > 2^{5}$ et $h_{n+1}^\circ \ge 24a_n h_n^\circ$, on peut construire un vecteur $\theta$ appartennant \`a $\mathbb{T}^2$ de mani\`ere que
\medskip

\noindent
(i) Si $$\sum_{n \ge 0} a_n^{-\frac{1}{3}} < \infty,$$ alors $T_\theta$ ne poss\`ede pas la propri\'et\'e du logarithme.
\medskip

\noindent
(ii) Si la suite $(a_n)$ admet une sous-suite born\'ee alors $\theta \in \Theta_l^2(0).$
\medskip

\noindent
(iii) Soit $\tau \ge 0$. On a $\theta \in  \Omega_s^2(\tau)$ si et seulement si $$ \inf a_{n+1}^{-1}a_n^{\frac{1+\tau}{2}}(h_n^\circ)^{1+ \tau}(h_{n+1}^\circ)^{-1} > 0.$$
%(iii) Soit $\tau \ge 0$. Si on a $$ \liminf a_{n+1}^{-1}a_n^{\frac{1+\tau}{2}}(h_n^\circ)^{1+ \tau}(h_{n+1}^\circ)^{-1} = 0,$$  alors  $\theta \notin  \Omega_s^2(\tau).$ Si on a $$\inf a_{n+1}^{-1}a_n^{\frac{1+\tau}{2}}(h_n^\circ)^{1+ \tau}(h_{n+1}^\circ)^{-1} >0,$$ alors  $\theta \in \Omega_s^2(\tau).$
%$h_{n+1}^\circ \le c a_{n+1}^{-1}a_{n}^{\frac{1+\tau}{2}} (h_n^\circ)^{1+\tau}$

\end{thm}
%$h_{n+1}^\circ =24 a_n h_n^\circ$
%c\frac{h_n^{\tau}}{a_{n+1}a_{n}^{\frac{1-\tau}{2}}}

Avant de faire la construction, on montre que l'on peut en d\'eduire le th\'eor\`eme \ref{thm4}. D'apr\`es la premi\`ere assertion du th\'eor\`eme, si on choisit $a_n = n^4$ (pour $n \ge 3$) alors $T_\theta$ ne poss\`ede pas la propri\'et\'e du logarithme. On choisit \'egalement $h_{n+1}^\circ =24 a_n h_n^\circ$. Alors, quelque soit $\tau>0$, $$ \inf a_{n+1}^{-1}a_n^{\frac{1+\tau}{2}}(h_n^\circ)^{1+ \tau}(h_{n+1}^\circ)^{-1} = \inf \frac{1}{24}(h_n^\circ)^\tau (n+1)^{-4} n^{2(\tau-1)} > 0,$$ car, $(h_n^\circ)$ a une croissance au moins exponentielle. Donc $\theta \in \Omega_s^2(\tau)$ pour tout $\tau > 0$. 
%Cela implique le \ref{thm4} $(i)$.

On choisit, maintenant, la suite $(a_n)$ born\'ee, donc $\theta \in \Theta_l^2(0)$ et $T_\theta$ poss\`ede la propri\'et\'e du logarithme. On veut de plus que $\theta$ n'appartienne \`a aucun $\Omega_s^2(\tau)$ pour $\tau >0$, soit d'apr\`es $(iii)$ que $ \inf a_{n+1}^{-1}a_n^{\frac{1+\tau}{2}}(h_n^\circ)^{1+ \tau}(h_{n+1}^\circ)^{-1} = 0$ pour tout $\tau >0$.  Comme la suite $(h_{n}^\circ)$ est soumise \`a la seule condition que $h_{n+1}^\circ \ge 24a_n h_n^\circ$, il suffit de la construire par r\'ecurrence en choisissant $h_{n+1}^\circ$ suffisamment grand devant $h_{n}^\circ$.
%Cela implique le \ref{thm4} $(ii)$.

%comme $(h_n^\circ)$ a une croissance au moins exponentielle, les suites $(h_n^\circ)$ et $(a_n)$ v\'erifient, quelque soit $\tau$ r\'eel, les conditions $(ii)$ et $(iii)$ du th\'eor\`eme \ref{thm7} ; cela implique le \ref{thm4} $(ii)$.

%D'apr\`es le $(i)$ du th\'eor\`eme $T_\theta$ ne poss\`ede pas la propri\'et\'e du logaritme,

%On chosit la suite $(a_n)$ born\'ee. Mais le choix de $h_{n+1}^\circ$ est ind\'ependant de $h_n$ et $a_n$. Si on construit par r\'ecurrence $h_{n+1}$ arbitrairement grand devant $a_n$, alors, comme $(h_n^\circ)$ a une croissance au moins exponentielle, $h_{n+1}^\circ$ est arbitrairement grand devant $q_n^\circ$. Ce qui implique aussi que $h_{n+1}$ est arbitrairement grand devant $q_n$   Il existe donc une $(h_n)$ et une suite $(a_n)$ qui respectivement v\'erifie la conditions du th\'eor\`eme 7 $(iii)$ et soit born\'ee.

%Si  $T_\theta$ poss\`ede la propri\'et\'e du logaritme et

%On se propose de construire une suite de droite $(\Delta_n)_{n \in \mathbb{N}}$ et une suite de points $(P_n)_{n \in \mathbb{N}}$, les droites $\Delta_n$ et $\Delta_{n+1}$ se coupant en $P_n$, de telle mani\`ere que les points $P_n$ convergent vers un point $\theta$, les suites $(\Delta_n)_{n \in \mathbb{N}}$ et $(P_n)_{n \in \mathbb{N}}$ nous donnant les meilleures approximations de $\theta$. Ce qui nous permettra de contr\^oler les qualit\'es diophantiennes de $\theta$ simultan\'ee et lin\'eaire.
%
%4.1
\subsection{Pr\'eliminaires}
La premi\`ere id\'ee est de se placer dans l'espace projectif. On consid\`ere les points $P = (x,y,z) \in \mathbb{Z}^3$, avec $x,y,z$ premiers entre eux, auxquels on associe les points $\widetilde{P}=\left( \frac{x}{z},\frac{y}{z} \right)$ dans $\mathbb{R}^2$ lorsque $z \ne 0$ et les triplets $\Delta=(r,s,t)$ de $\mathbb{Z}^3$, avec $r,s,t$ premiers entre eux associ\'es aux droites $ \widetilde{\Delta} : rx+sy+t=0$ dans $\mathbb{R}^2$.
\medskip

On notera $P \wedge P'$ le produit ext\'erieur dans $\mathbb{R}^3$ et $|P| = \max(|x|,|y|,|z|)$. Pour cette norme on a $|P \wedge P'| \le 2|P||P'|$.

Si $\widetilde{P}$, $\widetilde{P}'$ appartiennent \`a la droite $\widetilde{\Delta}$, ce qui correspond au fait que les points $P$ et $P'$ appartiennent au r\'eseau $\Gamma = \{ P \in \mathbb{Z}^3 | \langle \Delta, P \rangle = 0 \}$, alors $P \wedge P' = k \Delta$ avec $k$ entier et de plus $P$ et $P'$ engendrent le r\'eseau $\Gamma$ si et seulement si $P \wedge P' = \pm \Delta$.

On utilisera \'egalement la version suivante d'un lemme bien connu en norme euclidienne,
\begin{lem}
Soit $\Delta$ appartenant \`a $\mathbb{Z}^3$, le r\'eseau $\Gamma = \{ P \in \mathbb{Z}^3 | \langle \Delta, P \rangle = 0 \}$ et soit $P$ un \'el\'ement du r\'eseau, dont les coordonn\'ees sont premi\`eres entre elles. Alors, il existe $P'$ tel que $P$ et $P'$ engendrent $\Gamma$, avec $|P'| \le 2 \max \left( |P|, \frac{|\Delta|}{|P|}\right).$
\end{lem}
%
%4.2
\subsection{Construction de $(\Delta_n)_{n \in \mathbb{N}}$ et $(P_n)_{n \in \mathbb{N}}$ }

Nous allons construire $\theta$ comme limite d'une suite de points $(\widetilde{P}_n)$ tout en construisant une suite d'approximations lin\'eaires $(\widetilde{\Delta}_n)$.

Les normes respectives $h_n$ et $q_n$ de $\Delta_n$ et $P_n$ seront prescrites \`a une constante pr\`es. On demande de plus que le point $\widetilde{P}_n$ soit l'intersection de $\widetilde{\Delta}_n$ et de $ \widetilde{\Delta}_{n+1}$.
\medskip

Soient $(a_n)_{n \in \mathbb{N}}$ et $(h_n^\circ)_{n \in \mathbb{N}}$ deux suites d'entiers strictement positifs, avec $a_n > 2^{5}$ et $h_{n+1}^\circ \ge 24a_n h_n^\circ$. On note $q_n^\circ = a_n {h_n^\circ}^2$.
\medskip

Nous commen\c{c}ons la r\'ecurrence avec $\Delta_0= (h_0^\circ,-1,0)$ et $P_0 = (1,h_0^\circ, q_0^\circ)$, qui sont associ\'es respectivement \`a la droite $\widetilde{\Delta}_0 : h_0^\circ x-y=0$ et le point $\widetilde{P}_0 \left(\frac{1}{q_0^\circ} ,\frac{h_0^\circ}{q_0^\circ} \right) $. Alors $ q_0 = q_0^\circ$, $ h_0 = h_0^\circ$ et le point $\widetilde{P}_0$ appartient \`a la droite $\widetilde{\Delta}_0$. Supposons que soient construits $\Delta_n$ et $P_n$, \`a coordonn\'ees premi\`eres entre elles, tels que $\widetilde{P}_n$ appartient \`a $\widetilde{\Delta}_n$, et v\'erifiant les in\'egalit\'es $$\frac{1}{2}h_n^\circ \le | \Delta_n| = h_n \le 2h_n^\circ$$ et $$\frac{1}{2}q_n^\circ \le | P_n| = q_n \le 2q_n^\circ.$$

On construit d'abord $\Delta_{n+1}$. D'apr\`es le lemme $4.1$, il existe un point $\Delta'$ engendrant avec $\Delta_n$ le r\'eseau $\{ \Delta \in \mathbb{Z}^3 | \langle P_n, \Delta \rangle=0 \}$ et v\'erifiant $ | \Delta' | \le 2 \max(h_n,q_n h_n^{-1})$. On a $q_n h_n^{-1} \le 4a_nh_n^\circ$, donc $ | \Delta' | \le \frac{1}{3}h_{n+1}^\circ$. Quitte \`a changer de signe, on suppose $\Delta_n \wedge \Delta' = P_n$.

On pose
$$\Delta_{n+1} = \left[\frac{h_{n+1}^\circ}{h_n}\right]\Delta_n+\Delta' $$
et $h_{n+1} = |\Delta_{n+1}|$. On a encore $\Delta_n \wedge \Delta_{n+1} = P_n$, ce qui implique que les coordonn\'ees de $\Delta_{n+1}$ sont premi\`eres entre elles et que $\widetilde{P}_{n}$ appartient \`a $\widetilde{\Delta}_{n+1}$. De plus
$$ \left[\frac{h_{n+1}^\circ}{h_n}\right]h_n-|\Delta'| \le h_{n+1} \le \left[\frac{h_{n+1}^\circ}{h_n}\right]h_n+|\Delta'|,$$
$$ h_{n+1}^\circ-h_n -|\Delta'| \le h_{n+1} \le h_{n+1}^\circ +|\Delta'|.$$
D'o\`u
$$ \frac{1}{2}h_{n+1}^\circ \le h_{n+1} \le 2h_{n+1}^\circ.$$

On construit $P_{n+1}$ de la m\^eme mani\`ere ; le lemme $4.1$ donne l'existence d'un $P'$ engendrant avec $P_n$ le r\'eseau $\{ P \in \mathbb{Z}^3 | \langle \Delta_{n+1} ,P \rangle =0\}$ et v\'erifiant $ | P' | \le 2 \max(q_n,h_n q_n^{-1})$. Donc $ | P' | \le \frac{1}{3}q_{n+1}^\circ$. On pose alors $P_{n+1} = [\frac{q_{n+1}^\circ}{q_n}]P_n +P'$ et $q_{n+1}= |P_{n+1}|$. Alors le point $\widetilde{P}_{n+1}$ appartient \`a $\widetilde{\Delta}_{n+1}$ et on v\'erifie
$$ \frac{1}{2}q_{n+1}^\circ \le q_{n+1} \le 2q_{n+1}^\circ.$$
Quitte \`a changer de signe on supposera $P_n \wedge P_{n+1} = \Delta_{n+1}$.

%4.3
%
%
\subsection{Convergence de la suite $(\widetilde{P}_n)$ et d\'efinition de $\theta$}
Suivant la d\'emonstration de \cite{Laurent1}, on consid\`ere une \textit{``distance'' projective} entre  deux points $\widetilde{P}$ et $\widetilde{P'}$ de $\mathbb{R}^2$
$$d(\widetilde{P},\widetilde{P}')=\frac{|P \wedge P'|}{|P| |P'|}$$
et qui v\'erifie l'in\'egalit\'e triangulaire suivante pour trois points $\widetilde{P}$, $\widetilde{P}'$, $\widetilde{P}''$
$$ d(\widetilde{P},\widetilde{P}'') \le d(\widetilde{P},\widetilde{P}')+2d(\widetilde{P}',\widetilde{P}'').$$

Cette ``distance'' co\"incide avec la distance dans $\mathbb{R}^2$ pour deux points assez proches de l'origine. En effet, si $P= (x,y,z)$ et $P'=(x',y',z')$ v\'erifient $|z| \ge 2\max(|x|,|y|)$ et $|z'| \ge 2\max(|x'|,|y'|)$, un calcul imm\'ediat montre que $$d(\widetilde{P},\widetilde{P}')=|\widetilde{P}-\widetilde{P}'|.$$
Remarquons que $|z| \ge 2\max(|x|,|y|)$ si et seulement si $d(0,\widetilde{P}) \le \frac{1}{2}$, puisque $d(0,\widetilde{P}) = \frac{|(-y,x,0)|}{|P|}$.

Cela va nous permettre d'\'evaluer les distances $ |\widetilde{P}_n-\widetilde{P}_{n+1}|$. On a pour tout $n$ entier,
$$d(\widetilde{P}_n,\widetilde{P}_{n+1}) = \frac{|P_n \wedge P_{n+1}|}{|P_{n}| |P_{n+1}|}= \frac{|\Delta_{n+1}|}{|P_{n}| |P_{n+1}|}=\frac{h_{n+1}}{q_n q_{n+1}} \; \cdot$$
Donc
$$ \frac{1}{8} \frac{1}{a_{n+1}a_n h_{n}^\circ {}^2 h_{n+1}^\circ} \le \frac{1}{8}\frac{h_{n+1}^\circ}{q_n^\circ q_{n+1}^\circ} \le d(\widetilde{P}_n,\widetilde{P}_{n+1}) \le 8 \frac{h_{n+1}^\circ}{q_n^\circ q_{n+1}^\circ} \le \frac{8}{a_{n+1}a_n h_{n}^\circ {}^2 h_{n+1}^\circ}.$$
Compte tenu de l'encadrement pr\'ec\'edent et des conditions portant sur $(a_n)$ et $(h_n^\circ)$, on obtient

%\begin{equation*}
%d(\widetilde{P}_n,\widetilde{P}_{n+1}) \le \left (q_n^\circ 2^3 3^2 a_{n+1}a_n a_{n-1} h_{n-1}^\circ \right )^{-1}.
%\end{equation*}
%
%et
%\begin{equation*}
%d(\widetilde{P}_{n-1},\widetilde{P}_{n}) \ge 3\left (q_n^\circ a_{n-1} h_{n-1}^\circ \right)^{-1}.
%\end{equation*}
%Comme $a_n>2^5$, 
\begin{equation}
\label{inegpartie4large}
d(\widetilde{P}_n,\widetilde{P}_{n+1}) \le \frac{1}{2^{18}3^{3} }d(\widetilde{P}_{n-1},\widetilde{P}_{n}).
\end{equation}
%En particulier, on a $$ $$

Pour $n=0$, on a $d(0,\widetilde{P}_0) = \frac{h_0^\circ}{q_0^\circ} = \frac{1}{a_0h_0^\circ} < \frac{1}{32}$. Et en \'evaluant grossi\`erement, $d(\widetilde{P_{i}},\widetilde{P}_{i+1}) \le \frac{1}{32^{i+1}}$, d'o\`u,
$$d(0,\widetilde{P}_n) \le d(0,\widetilde{P}_0) + 2\sum_{i=0}^{n-1} d(\widetilde{P}_i,\widetilde{P}_{i+1}) < \frac{1}{32} + 2\sum_{i=0}^{n-1} \frac{1}{32^{i+1}} < \frac{1}{8}.$$
En particulier, $d(0,\widetilde{P}_n) < \frac{1}{2}$, pour tout $n$ et
$$ |\widetilde{P}_n-\widetilde{P}_{n+1}| = d(\widetilde{P}_n,\widetilde{P}_{n+1}) .$$
%= \frac{h_{n+1}}{q_n q_{n+1}}

La d\'ecroissance de $|\widetilde{P}_n-\widetilde{P}_{n+1}|$ \'etant au moins g\'eom\'etrique, la suite $\widetilde{P}_n$ est de Cauchy et on note $\theta$ sa limite. On a encore 
\begin{equation}
\label{inegtheta}
|\theta|\le \frac{1}{8} < \frac{1}{4}.
\end{equation}

De (\ref{inegpartie4large}), il r\'esulte aussi $\frac{1}{2}|\widetilde{P}_n-\widetilde{P}_{n+1}| \le |\widetilde{P}_n-\theta| \le \frac{3}{2} |\widetilde{P}_n-\widetilde{P}_{n+1}|,$ soit encore

\begin{equation}
\label{ineg5partie4}
\frac{1}{2} \frac{h_{n+1}}{q_n q_{n+1}} \le |\widetilde{P_n}-\theta| \le \frac{3}{2} \frac{h_{n+1}}{q_n q_{n+1}}.
\end{equation}
%\frac{1}{2} \frac{h_{n+1}}{q_n q_{n+1}} \le
%
%
\subsection{Meilleures approximations diophantiennes de $\theta$}
%4.4.1
\subsubsection{Approximations simultan\'ees}
Nous montrons ici que les meilleures approximations simultan\'ees de $\theta=(\theta',\theta'')$ sont les $|P_n|=q_n$ \`a l'exception possible des premiers termes et de certains termes $(q_{n+1}-q_n)$ :%

\begin{lem}
Si $q$ est un entier avec $0 < q < q_{n+1}$ et $q \ne q_n$, $q \ne q_{n+1}-q_n$, alors
$$\| q \theta \| > \| q_n \theta \| . $$
\end{lem}

%Si $\widetilde{\Delta}$ passe par $\widetilde{P_n}$ avec $\widetilde{\Delta} \ne \widetilde{\Delta_n}$, alors $ | \Delta | >\frac{1}{16}a_n h_n > | \Delta_n|.$

%Les conditions se r\'e\'ecrivent, $\langle \Delta, P_n \rangle=0$ et $\Delta_n$ non colin\'eaire \`a $\Delta$.
%Comme $\langle \Delta_n, P_n \rangle=0$, on a $ (\Delta_n \wedge \Delta) \wedge P_n =0$.

%On a donc pour un entier $k$, $ \Delta_n \wedge \Delta = kP_n$ ; or $ \Delta_n \wedge \Delta \ne 0$, donc $k$ est non nul.
%On en d\'eduit donc ,
%$$| \Delta | \ge \frac{|P_n|}{2|\Delta_n|} = \frac{q_n}{2h_n} \ge \frac{1}{16}a_nh_n > |\Delta_n|.\qed$$
%

D'apr\`es (\ref{ineg5partie4}), on a $q_n| \theta - \widetilde{P}_n| \le \frac{1}{2},$ donc $ \| q_n \theta \| = q_n| \theta - \widetilde{P}_n |$ et
\begin{equation} \label{ineg6partie4}  \frac{h_{n+1}}{2q_{n+1}} \le \| q_n \theta \| \le \frac{3h_{n+1}}{2q_{n+1}} ,\end{equation}
%
%
%
% car $q_{n+1} \ge \frac{1}{2} q_{n+1}^\circ= \frac{1}{2} a_{n+1}((h_{n+1}^\circ)^2 \ge \frac{1}{8}a_{n+1}h_{n+1}^2.$

Soit $q$ un entier satisfaisant les condition du lemme. Il existe $x$ et $y$ entiers tels que  et $\| q \theta \| = \max(|x-q\theta'|, |y-q\theta''|)$. Soit $P =(x,y,q)$. D'apr\`es (\ref{inegtheta}), $|x| < \frac{q}{4} +\frac{1}{2} < q$ et de m\^eme $|y| < q$. On a donc $|P|=q$. Si $|\widetilde{P}|>\frac{1}{2}$ alors on a $| \widetilde{P} - \theta | > \frac{1}{4}$, or d'apr\`es (\ref{ineg6partie4}), $ \|q_n\theta \| < \frac{1}{4}$. Nous supposerons donc que $|\widetilde{P}| \le \frac{1}{2}$,

On d\'ecompose
$$ \| q \theta \| = q | \widetilde{P} - \theta | \ge q \left ( | \widetilde{P} - \widetilde{P}_{n+1}| - | \widetilde{P}_{n+1} - \theta | \right)$$

D'apr\`es (\ref{ineg5partie4}),
$ q | \widetilde{P}_{n+1} - \theta | \le \frac{3}{2} \frac{h_{n+2}}{q_{n+2}},$ d'o\`u
%Rappelons que l'on a, $\frac{1}{8}a_nh_n2 \le q_n \le 8a_nh_n2$ et $h_{n+1} \ge 6 a_nh_n$.
%$$|P|d(\widetilde{P},\theta) \ge \frac{|P \wedge P_{n+1}|}{|P_{n+1}|} - \frac{25 |P|}{a_{n+2}h_{n+2}q_{n+1}}.$$
$$ q | \widetilde{P}_{n+1} - \theta | \le \frac{6}{a_{n+2}h_{n+2}^\circ} \le \frac{1}{4 \, a_{n+2}a_{n+1}h_{n+1}^\circ} \le \frac{h_{n+1}^\circ}{4 \, a_{n+2}q_{n+1}^\circ} \le \frac{ h_{n+1}}{a_{n+2}q_{n+1}},$$
donc
$$ q | \widetilde{P}_{n+1} - \theta | < \frac{1}{2} \frac{h_{n+1}}{q_{n+1}}.$$
Pour conclure la d\'emonstration du lemme, il sufit maintenant de montrer que $$q| \widetilde{P} - \widetilde{P}_{n+1}| \ge 2 \, \frac{h_{n+1}}{ q_{n+1}} .$$
Rappelons que l'on suppose $|\widetilde{P}| \le \frac{1}{2}$,  on a donc
\begin{equation}
\label{ineg9partie4}
 | \widetilde{P} - \widetilde{P}_{n+1}| = \frac{| P \wedge P_{n+1}|}{|P| |P_{n+1}|} = \frac{1}{q \, q_{n+1}}| P \wedge P_{n+1}| .
\end{equation}

Il reste \`a \'evaluer $| P \wedge P_{n+1}|$.

Premi\`ere possibilit\'e, $\widetilde{P}$ n'appartient pas \`a $\widetilde{\Delta}_{n+1}$, auquel cas la droite joignant $\widetilde{P}$ et $\widetilde{P}_{n+1}$ passe par $\widetilde{P}_{n+1}$ tout en n'\'etant pas colin\'eaire \`a $\widetilde{\Delta}_{n+1}$. Il existe alors un entier $l$ non nul tel que $$(P \wedge P_{n+1}) \wedge \Delta_{n+1} = lP_{n+1}.$$
D'o\`u
$$|P \wedge P_{n+1}| \ge \frac{|P_{n+1}|}{2|\Delta_{n+1}|} = \frac{q_{n+1}}{2h_{n+1}} \ge 2h_{n+1}.$$

Deuxi\`eme possibilit\'e, $\widetilde{P}$ appartient \`a $ \widetilde{\Delta}_{n+1}$, c'est \`a dire $P$ appartient au r\'eseau $\{ P \in \mathbb{Z}^3 | \langle \Delta_{n+1} ,P \rangle =0\}$. Les points $P_n$ et $P_{n+1}$ engendrant ce r\'eseau, il existe $k$ et $k'$ entiers tels que $P = kP_n + k' P_{n+1}$. On a
$$| P \wedge P_{n+1}| = |k| |\Delta_{n+1}|=|k|h_{n+1}.$$
Remarquons que $k$ est non nul car $ |P| < |P_{n+1}|$. Si $|k| = 1$, on a alors  $P =  k'P_{n+1} \pm P_n $ et $q =    k'q_{n+1} \pm q_n $. Comme  $ 0 < q < q_{n+1}$ on a $k' = 0$ ou $k' = 1$, mais comme $q \ne q_n$ et  $ q \ne  q_{n+1} - q_n $, ces cas sont \'egalement exclus. D'o\`u $|k| \ge 2$ et on a encore 
$$| P \wedge P_{n+1}|\ge 2h_{n+1}.$$

D'apr\`es (\ref{ineg9partie4}), $q| \widetilde{P} - \widetilde{P}_{n+1}| \ge 2 \, \frac{h_{n+1}}{ q_{n+1}},$ ce qui conclut la d\'emonstration du lemme.
%

%4.4.2
\subsubsection{Approximations lin\'eaires}
Notons ici $\bar{\theta}= (\theta', \theta'',1) \in \mathbb{R}^3$ et $ \Delta_n = (r_n, s_n, t_n ) $, alors $ | \langle \Delta_n, \bar{\theta} \rangle | = | r_n \theta' + s_n\theta'' +t_n |.$

Nous allons montrer que $(r_n, s_n)_{n \in \mathbb{N}}$ est la suite des meilleures approximations lin\'eaires de $\theta$, \`a l'exception possible des premiers termes et de certains $(r_{n+1},s_{n+1}) \pm (r_{n},s_{n})$.
\medskip

Donnons tout d'abord un encadrement de $ | \langle \Delta_n,\bar{\theta} \rangle | $.
\begin{lem}
\label{encadrementDeltan}Pour tout $n$ entier, on a 
$$\frac{3}{4}\frac{1}{q_{n+1}} \le | \langle \Delta_n,\bar{\theta} \rangle | \le \frac{5}{4} \frac{1}{q_{n+1}}.$$
\end{lem}

On d\'ecompose, $$ | \langle \Delta_n,\bar{\theta} \rangle | = \left| \langle \Delta_n,\frac{P_{n+1}}{|P_{n+1}|}\rangle + \langle \Delta_n, \bar{\theta} -\frac{P_{n+1}}{|P_{n+1}|} \rangle \right|.$$
Etudions $ \langle \Delta_n,P_{n+1} \rangle$, on a
$$\Delta_n \wedge (P_n \wedge P_{n+1}) = \Delta_n \wedge \Delta_{n+1} = P_n.$$
Mais aussi, d'apr\`es la formule classique sur le double produit vectoriel
$$\Delta_n \wedge (P_n \wedge P_{n+1}) = \langle \Delta_n,P_{n+1} \rangle P_n - \langle \Delta_n,P_{n} \rangle P_{n+1} = \langle \Delta_n,P_{n+1} \rangle P_n.$$ Donc, $ \langle \Delta_n,P_{n+1} \rangle= 1,$ et
$$\left| \langle \Delta_n,\frac{P_{n+1}}{|P_{n+1}|} \rangle \right| = \frac{1}{q_{n+1}}.$$
D'autre part d'apr\`es (\ref{ineg5partie4}),
$$\left| \langle \Delta_n, \bar{\theta}-\frac{P_{n+1}}{|P_{n+1}|} \rangle \right| \le 2 h_n | \theta- \widetilde{P}_{n+1} | \le \frac{3h_n h_{n+2}}{q_{n+1}q_{n+2}}  \le \frac{24h_n}{a_{n+2}q_{n+1}h_{n+2}} \le \frac{1}{4} \frac{1}{q_{n+1}}.$$
D'o\`u le r\'esultat.
\medskip
\medskip
\medskip
\medskip

En particulier,  $ | \langle \Delta_n, \bar{\theta} \rangle | < \frac{1}{2} $, donc
$$ | \langle \Delta_n,\bar{\theta} \rangle | = | r_n \theta' + s_n\theta'' +t_n | = \| \langle (r_n,s_n), \theta \rangle \|.$$
Pour montrer qu'il s'agit des meilleures approximations, nous allons comparer  $\Delta_n$ \`a d'autres vecteurs $\Delta$ de $\mathbb{Z}^3$.  Rappelons que $$|\theta|=\max(|\theta'|,|\theta''|) < \frac{1}{4}.$$ 
Comme $| \langle \Delta_n,\bar{\theta} \rangle | < \frac{1}{2}$, on a $|t_n| < \frac{1}{2} +\frac{1}{2}|(r_n,s_n) | <|(r_n,s_n)| $ et donc $|(r_n,s_n)| =h_n$. Soient $(r,s)$ un vecteur du tore tel que $0 < |(r,s)| < h_{n+1}$ et $t$ entier tel que $ \| \langle (r,s),\theta) \rangle \| = |r \theta' + s\theta''+t|$. En notant  $\Delta = (r,s,t)$, $ \langle \Delta,\bar{\theta} \rangle  =  \langle (r,s),\theta) \rangle$ et de m\^eme $|\Delta | = |(r,s)| <h_{n+1}$.

Pour v\'erifier que $(r_n,s_n)$ est la suite des meilleures approximations lin\'eaires de $\theta$ \`a l'exception possible des premiers termes et de certains termes $\pm (r_{n+1},s_{n+1}) \pm (r_{n},s_{n})$, il nous suffit donc de montrer le lemme suivant.
\begin{lem}
Soit $\Delta = (r,s,t)$ un vecteur de $\mathbb{Z}^3$ tel que $| \Delta | < h_{n+1}$, $ \Delta \ne \pm \Delta_n$ et $ \Delta \ne \pm \Delta_{n+1} \pm \Delta_n$. Alors $$| \langle \Delta, \bar{\theta} \rangle | > | \langle \Delta_n, \bar{\theta} \rangle | .$$
\end{lem}

On d\'ecompose de la m\^eme mani\`ere que dans le lemme $4.3$ et on obtient
$$\left| \langle \Delta, \bar{\theta}-\frac{P_{n+1}}{|P_{n+1}|}\rangle \right| \le \frac{24h_{n+1}}{a_{n+2}q_{n+1}h_{n+2}} \le \frac{1}{4} \frac{1}{q_{n+1}}.$$
De plus ,
$$\left| \langle \Delta, \frac{P_{n+1}}{|P_{n+1}|} \rangle \right| = \frac{1}{q_{n+1}} |\langle\Delta,P_{n+1}\rangle|.$$
D'o\`u
$$\left| \langle \Delta, \bar{\theta} \rangle \right| \ge \frac{1}{q_{n+1}} \left ( |\langle\Delta,P_{n+1}\rangle| - \frac{1}{4} \right ). $$

D'apr\`es le lemme \ref{encadrementDeltan}, il suffit donc de montrer que $$\left| \langle \Delta, P_{n+1} \rangle \right| \ge 2.$$

Tout d'abord remarquons que $\langle \Delta, P_{n+1} \rangle$ est non nul. Dans le cas contraire, $\widetilde{P}_{n+1}$ appartiendrait \`a $\widetilde{\Delta}$, qui est non colin\'eaire \`a $\widetilde{\Delta}_{n+1}$ puisque $|\Delta| < h_{n+1}$. Alors $\Delta \wedge \Delta_{n+1} = k P_{n+1}$, avec $k$ un entier non nul, et donc
$$|\Delta| \ge \frac{q_{n+1}}{2h_{n+1}} \ge \frac{a_n h_{n+1}}{16} \ge 2h_{n+1},$$
ce qui est contraire \`a l'hypoth\`ese.

Supposons maintenant que, $\langle \Delta,P_{n+1}\rangle = \pm 1$. Quitte \`a changer $\Delta$ en $-\Delta$, on peut supposer $\langle \Delta,P_{n+1}\rangle = \langle \Delta_n,P_{n+1}\rangle$, d'o\`u $\langle \Delta-\Delta_n,P_{n+1}\rangle = 0$.

Alors la droite correspondant au triplet $\Delta - \Delta_n $, passe par $\widetilde{P}_{n+1}$. Si cette droite n'est pas colin\'eaire \`a $ \widetilde{\Delta}_{n+1}$, on en d\'eduit de m\^eme
$$|\Delta-\Delta_n| \ge \frac{q_{n+1}}{2h_{n+1}} \ge 2h_{n+1},$$
ce qui est de nouveau impossible, puisque $h_n < h_{n+1}$. On a donc $\Delta = \Delta_n + k \Delta_{n+1} $, avec $k$ entier.

Si $k=0$ ou $k= \pm 1$, on a respectivement  $\Delta = \Delta_n$ ou $\Delta = \pm \Delta_{n+1} + \Delta_n$, ces deux cas ont \'et\'e exclus. De plus $|k| \ge 2$ impliquerait $|\Delta| > h_{n+1}$, ce qui est en contradiction avec notre hypoth\`ese.

Donc $\left| \langle \Delta, P_{n+1} \rangle \right| \le 1$ est impossible, ce qui conclut la d\'emonstration du lemme.

%.4.5
\subsection{D\'emonstration du th\'eor\`eme $7$ }
Nous avons donn\'e une construction d'un vecteur $\theta$ permettant d'en contr\^oler les meilleures approximations, selon le choix de deux suites $(a_n)$ et $(h_n^\circ)$. Il ne nous reste plus qu'\`a pr\'eciser les types diophantiens selon  $(a_n)$ et $(h_n^\circ)$ et v\'erifier si $\theta$ poss\`ede ou non la propri\'et\'e du logarithme.

D'apr\`es le lemme \ref{encadrementDeltan}, on a
$$\frac{3}{4} h_{n+1}^2\frac{1}{q_{n+1}}  \le |\Delta_{n+1}|^2 |\langle \Delta_n, \bar{\theta} \rangle | \le \frac{5}{4} h_{n+1}^2\frac{1}{q_{n+1}},$$
$$ \frac{3}{32}\frac{1}{a_{n+1}} \le |\Delta_{n+1}|^2 |\langle \Delta_n, \bar{\theta} \rangle | \le \frac{10}{a_{n+1}}.$$
En dimension 2, en utilisant les $(\Delta_n)$, la condition du th\'eor\`eme $5$ s'\'ecrit $$ \sum (|\Delta_{n+1}|^2 |\langle \Delta_n, \bar{\theta} \rangle |)^{\frac{1}{3}} < \infty.$$ Donc si la suite $(a_n)$ v\'erifie $\sum a_n^{-\frac{1}{3}} < \infty$, alors $T_\theta$ ne poss\`ede pas la propri\'et\'e du logarithme. Par contre, si la suite $(a_n)$ a une sous-suite born\'ee, alors d'apr\`es l'in\'egalit\'e de gauche, il existe une constante $c$ strictement positive telle que pour tout $n$ $$|\Delta_{n+1}|^2 \eps_l(|\Delta_n|) \ge c,$$
d'o\`u $\theta \in \Theta_l^2(0)$.
%Comme la suite $(\Delta_n)$ correspond, au moins, \`a une sous-suite de la suite des meilleures approximations, s
%$$\limsup_{h \rightarrow +\infty} h^{d(1+\tau)} \ \eps_l(h) >0,$$ .$
%
\par\medskip\noindent
Nous avons montr\'e les assertions $(i)$ et $(ii)$ du th\'eor\`eme.
%Nous aurions aussi pu utiliser les approximations simultan\'ees, l'in\'equation (\ref{ineg6partie4}) et la proposition \ref{prop32} nous donneraient alors que $T_\theta$ ne poss\`ede pas la propri\'et\'e du logarithme si $\sum a_n^{-\frac{1}{6}} < \infty$. On montrerait \'egalement que si $(a_n)$ a une sous suite born\'ee, alors $\theta \in \Theta_s^2(0)$.
\medskip

On cherche maintenant les conditions sur $(a_n)$ et $(h_n^\circ)$ pour que $\theta$ appartiennent ou non \`a $\Omega_s^2(\tau)$ pour un $\tau \ge 0$ donn\'e. D'apr\`es l'encadrement (\ref{ineg6partie4}), on a
$$\frac{1}{8} \frac{1}{a_{n+1}h_{n+1}^\circ} \le \| q_n \theta \| \le \frac{6}{a_{n+1}h_{n+1}^\circ}, $$
$$\frac{1}{8}\frac{q_n^{\frac{1+\tau}{2}}}{a_{n+1}h_{n+1}^\circ} \le q_n^{\frac{1+\tau}{2}} \| q_n \theta \| \le 6\frac{q_n^{\frac{1+\tau}{2}}}{a_{n+1}h_{n+1}^\circ} \; , $$

\begin{equation}
\label{inega7partie4}
c_0 \frac{a_n^{\frac{1+\tau}{2}} (h_n^\circ)^{1+ \tau}}{a_{n+1}h_{n+1}^\circ} \le q_n^{\frac{1+\tau}{2}} \| q_n \theta \| \le c_1\frac{a_n^{\frac{1+\tau}{2}} (h_n^\circ)^{1+ \tau}}{a_{n+1}h_{n+1}^\circ} \; ,
\end{equation}
o\`u $c_0$ et $c_1$ sont deux constantes strictement positives.
\medskip

Pour les \'eventuelles meilleures approximations de la forme $q_{n+1}- q_n$, remarquons que 
$$(q_{n+1}- q_n)^{\frac{1+\tau}{2}}\eps_s(q_{n+1}- q_n) \ge 2^{-\frac{1+\tau}{2}} q_{n+1}^{\frac{1+\tau}{2}} \| q_{n+1} \theta \|,$$
puisque $\eps_s(q_{n+1}- q_n) \ge \eps_s(q_{n+1})$ et $q_{n+1}- q_n \ge \frac{1}{2}q_{n+1}$.

Donc $\theta \in \Omega_s^2(\tau)$ si et seulement si $\inf q_{n}^{\frac{1+\tau}{2}} \| q_{n+1} \theta \| > 0$ et d'apr\`es (\ref{inega7partie4}), cela est vrai  si et seulement si 
$$ \inf a_n^{\frac{1+\tau}{2}}(h_n^\circ)^{1+ \tau}(a_{n+1}h_{n+1}^\circ)^{-1} > 0.$$
Ce qui conclut la d\'emonstration du th\'eor\`eme.
%$\liminf q^{\frac{1+\tau}{d}}\eps_s(q) >0$

%%%%%%%%%%%%%%%%%%%

%
\medskip
\medskip

\noindent
LAGA, UNIVERSIT\'E PARIS 13, CNRS UMR 7539.
\medskip

\noindent
\textit{E-mail address} : mussat@math.univ-paris13.fr

\end{document}